\documentclass[11pt]{article}
\usepackage{amsfonts, amssymb, graphicx, a4}
\usepackage{amsmath}
 \usepackage{color}

\newtheorem{theorem}{Theorem}[section]
\newtheorem{lemma}[theorem]{Lemma}
\newtheorem{proposition}[theorem]{Proposition}

\newtheorem{definition}[theorem]{Definition}
\newtheorem{corollary}[theorem]{Corollary}
\newtheorem{remark}[theorem]{Remark}
\newtheorem{remarks}[theorem]{Remarks}

 \newenvironment{proof}{\noindent {\bf Proof.\,}
 }{\hspace*{\fill}$\qed$\medskip}

 \newenvironment{rems}{\begin{remarks}\rm}{\end{remarks}}

\newcommand{\floor}[1]{\lfloor#1\rfloor}

\def \P {{\mathbb P}}
\def \E {{\mathbb E}}

\newcommand{\cX}{{\cal X}}

\def \a {{\alpha}}
\def \b {{\beta}}
\def \e {{\varepsilon}}
\def \D {{\Delta}}
\def \r {{\rho}}
\def \l {{\lambda}}

\def \z {{\zeta}}

\def \t {{\tau}}

\def \d {{\delta}}

\def\eqref#1{(\ref{#1})}

\newcommand{\be}[1]{\begin{equation}\label{#1}}
\newcommand{\ee}{\end{equation}}

\newcommand{\bl}[1]{\begin{lemma}\label{#1}}
\newcommand{\el}{\end{lemma}}

\newcommand{\br}[1]{\begin{remark}\label{#1}}
\newcommand{\er}{\end{remark}}

\newcommand{\bt}[1]{\begin{theorem}\label{#1}}
\newcommand{\et}{\end{theorem}}

\newcommand{\bd}[1]{\begin{definition}\label{#1}}
\newcommand{\ed}{\end{definition}}

\newcommand{\bcl}[1]{\begin{claim}\label{#1}}
\newcommand{\ecl}{\end{claim}}

\newcommand{\bp}[1]{\begin{proposition}\label{#1}}
\newcommand{\ep}{\end{proposition}}

\newcommand{\bc}[1]{\begin{corollary}\label{#1}}
\newcommand{\ec}{\end{corollary}}

\newcommand{\bi}{\begin{itemize}}
\newcommand{\ei}{\end{itemize}}

\newcommand{\ben}{\begin{enumerate}}
\newcommand{\een}{\end{enumerate}}
 \newcommand{\bpr}{\begin{proof}}
 \newcommand{\epr}{\end{proof}}

\def\tc{ : \ }
\def \qed {{\square\hfill}}

\def \P {{\mathbb P}}
\def \E {{\mathbb E}}

\def\1{\rlap{\mbox{\small\rm 1}}\kern.15em 1}

\def\build#1_#2^#3{\mathrel{\mathop{\kern 0pt#1}\limits_{#2}^{#3}}}
\def\tend#1#2#3{\build\hbox to 12mm{\rightarrowfill}_{#1\rightarrow
#2}^{#3}}

\def\tendn{\tend{n}{\infty}}{}
\def\converge#1#2#3{\build\hbox to
15mm{\rightarrowfill}_{\hbox{\scriptsize #3}}^{#1\rightarrow #2}}
\def\converg#1#2#3{\build\hbox to
15mm{\rightarrowfill}_{\hbox{\scriptsize #3}}^{#1\uparrow #2}}

\newcommand{\eproof}{\hspace*{\fill}$\qed$}

\newcommand{\empbf}[1]{{\bf \emph{#1}}}

\begin{document}
\author{
\renewcommand{\thefootnote}{\arabic{footnote}}
R.\ Fernandez\,\footnotemark[1]\;,
\renewcommand{\thefootnote}{\arabic{footnote}}
F.\ Manzo\,\footnotemark[2]\;,
\renewcommand{\thefootnote}{\arabic{footnote}}
F. R.\
Nardi\,\footnotemark[3]\;\textsuperscript{\,,}\footnotemark[4]\;, and
\renewcommand{\thefootnote}{\arabic{footnote}}
E.\ Scoppola\,\footnotemark[2]
}

\title{Asymptotically exponential hitting times and metastability:\\
 a pathwise approach without reversibility}

\footnotetext[1]{
Mathematics Department, Utrecht University, P.O.\ Box 80010, 3508 TA
Utrecht, The Netherlands
}
\footnotetext[2]{
Dipartimento di Matematica,
 Universit\`a di Roma Tre,
 Largo S.\ Leonardo Murialdo 1,
 00146 Rome, Italy
}
\footnotetext[3]{
Technische Universiteit Eindhoven, P.O.\ Box 513, 5600 MB Eindhoven,
The Netherlands
}
\footnotetext[4]{
EURANDOM, P.O.\ Box 513, 5600 MB Eindhoven, The Netherlands
}

\maketitle


%
%

\begin{abstract}
We study the hitting times of Markov processes  to target set $G$, starting from a reference configuration $x_0$ or its basin of attraction. The configuration $x_0$ can correspond to the bottom of a (meta)stable well, while the target $G$ could be either a set of saddle (exit) points of the well, or a set of further (meta)stable configurations. Three types of results are reported:
(1) A general theory is developed, based on the path-wise approach to metastability,
which has three important attributes.  First, it is general in that it does not assume reversibility of the process, does not focus only on hitting times to rare events and does not assume a particular starting measure.  Second, it relies only on the natural hypothesis that the mean hitting time to $G$ is asymptotically longer than the mean recurrence time to $x_0$ or $G$. Third, despite its mathematical simplicity, the approach yields precise and explicit bounds on the corrections to exponentiality.
(2) We compare and relate different metastability conditions proposed in the literature so to eliminate potential sources of confusion.  This is specially relevant for evolutions of infinite-volume systems, whose treatment depends on whether and how relevant parameters (temperature, fields) are adjusted.
(3) We introduce the notion of early asymptotic exponential behavior to control time scales asymptotically smaller than the mean-time scale. This control is particularly relevant for systems with unbounded state space where nucleations leading to exit from metastability can happen anywhere in the volume. We provide natural sufficient conditions on recurrence times for this early exponentiality to hold and show that it leads to estimations of probability density functions.


\end{abstract}

\eject


\section{Introduction}

Hitting times to rare sets, and related metastability
issues, have been studied
both in the framework of probability theory  and statistical mechanics.

A short review of first hitting results is given in sect. 1.1.
As far as metastability results are concerned, the story is much more involved due to the fact that metastability can be defined in different ways.

The phenomenon of metastability is given by the following scenario: (i) A system remains ``trapped'' for an abnormally long time in a state ---the metastable phase--- different from the eventual equilibrium state consistent with the thermodynamical potentials. (ii) Subsequently, the system undergoes a \emph{sudden} transition  from the metastable to the stable state at a \emph{random time}.

The mathematical study of this phenomenon has been a standing issue since the foundation of the field of rigorous statistical mechanics.  This long history resulted in a number of different approaches, based on non-equivalent assumptions.

The common feature of these mathematical descriptions is a stochastic framework involving a Markovian evolution and two disjoint distinguished sets of configurations, respectively representing metastable and stable states. Rigorously speaking, a statistical mechanical ``state'' corresponds to a probability measure.  The association to sets of configurations corresponds, therefore, to identifying supports of relevant measures.   This is a crucial step in which both probabilistic and physical information must be incorporated. Within such  framework, the central mathematical issue is the description of the first-exit trajectories leading from an initial metastable configuration to a final stable one.
 The exit path can be decomposed into two parts:
 \begin{itemize}
\item an \emph{escape path} ---taking the system to the boundary of the metastable set, which can be  thought of as a saddle or a bottleneck: a set with small equilibrium measure which is difficult for the system to reach. ---
\item and a \emph{downhill path} ---bringing the system into the next stable state.  ---
 \end{itemize}
Metastable behavior  occurs when the time spent in the downhill path is negligible  with respect to the
escape time.
In such a scenario, the overall exponential character of the time to reach stability is, therefore, purely due to the escape part of the trajectory.  This implies that the set of metastable states can be considered as a
single state in the sense that  the first escape  turns out to be
exponentially distributed as in the case of a single state.
Actually this exponential law can be considered as the main feature of metastability.

As noted above the first escape from the metastable set can be seen as the first hitting to a rare
set of saddle configurations   and this is the strong
relation between metastability and first hitting to rare events. However we have to note that
the important quantity in metastability is the decay time of the metastable phase so the first hitting
to the stable one. In order to reduce this problem to the first hitting to a rare set, the saddle set, one has
to individuate this boundary set by investigating in detail the state space.

In metastability literature, the main used tools are
 renormalization \cite{Scop93,Scop94}, cycle decomposition and large deviations \cite{CaCe,Ca,CaTr,CGOV,OS1,OS2,OV,MNOS,Trouve1}
and, lately, potential theoretic techniques \cite{BovMetaPTA,BEGK00,BEGK02,BEGK04,BL1,BG}.
Generally speaking, the focus is more on the exit path and on the mean exit time, while the results on the distribution of the escape time are usually asymptotical and not quantitative as in  \cite{AB1}.

At any rate, metastability involves a relaxation from an initial measure to a metastable state, from which the systems undergoes a final transition into stability.  To put this two-step process in evidence, the theory must apply to sufficiently general initial states.

This paper is based on the \emph{path-wise approach} to metastability developed in\cite{CGOV,OS1,OS2,OV,MNOS}.  This approach has led to a detailed description of exit paths in terms of relevant physical quantities, such as energy, entropy, temperature, etc.
In this paper we show how the same approach provides simple and effective estimations of the laws of hitting times for rather general, not necessarily reversible dynamics.
Moreover  we can consider hitting to more general goals, not necessarily rare sets,
so, in metastability
language, we do not need  to individuate the saddle configurations
in order to prove the exponentially of the decay time since we can directly consider as goal the
stable state.
This can be an  important point to apply the
theory in very complicated physical contexts where is too difficult a detailed control of the state space.

 \subsection{First hitting: known results}

The stochastic treatment was initially developed in the framework of  reliability theory, in which the reference states are called \emph{good} states and the escaped states are the \emph{bad} states.
The exponential character of good-to-bad transitions ---well known from quite some time--- is due to the existence of two different time scales:   Long times are needed to go from good to bad states, while the return to good states from anywhere ---except, perhaps, the bad states--- is much shorter.  As a result, a system in a good state can arrive to the bad state only through a large fluctuation that takes it all the way to the bad state.  Any intermediate fluctuation will be followed by an unavoidable return to the good states, where, by Markovianness, the process starts afresh independently of previous attempts.  The escape time is formed, hence, by a large number of independent returns to the good states  followed by a final successful excursion to badness that must happen without hesitations, in a much shorter time.  The exit time is, therefore, a geometric random variable with extremely small success probability; in the limit, exponentiality follows.

A good reference to this classical account is the short book \cite{Kei79} which, in fact, collects also the main tools subsequently used in the field: reversibility, spectral decomposition, capacity, complete monotonicity. Exponentiality of hitting times to rare events is analyzed, in particular, in Chapter 8 of this book, where regenerative processes are considered.

The tools given in \cite{Kei79} where exploited in \cite{B1,B2,B3,A,AB1,AB2} to  provide  sharp and explicit estimates
on the exponential behavior of hitting times with means larger than the relaxation time of the chain.  For future reference, let us review some results of these papers.

Let $X_t;\,  t\ge 0$ be an irreducible, finite-state, reversible Markov chain in continuous time, with
 transition rate matrix $Q$ and stationary distribution $\pi$.  The {relaxation time} of the chain is
 $R=1/\l_1$, where $\l_1$ is the smallest eigenvalue of $-Q$.  Let $\t_A$ denote the hitting time of a given subset $A$  of the state space and $\P_\pi$ the law of the process started at $\pi$.  Then,
 for all $t>0$ (Theorem 1 in  \cite{AB1}):
 \be{th1AB}
 \Bigl|\P_\pi\bigl(\t_A/\E_\pi\t_A>t\bigr)-e^{-t}\Bigr|\;\le\;
  {R/\E_\pi\t_A\over 1+R/\E_\pi\t_A}\;.
 \ee
 Moreover in the regime $R\ll t\ll \E_\pi\t_A$,  the distribution of $\t_A$ rescaled by its mean value, can be can be controlled with explicit bounds on its density function (Sect. 7 of \cite{AB1}).
 We note that these results are given in terms of quantitative estimates which are strict when $R/\E_\pi\t_A\ll 1$ and the starting
 distribution is the equilibrium measure.

Further insight is provided by the \emph{quasi-stationary distribution}
 $$
 \a\; :=\;\lim_{t\to\infty}\P_\pi(X_t\in \cdot\mid\t_A>t)\;.
 $$
This distribution is stationary for the conditional process
$$
\a\;=\;\P_\a(X_t\in \cdot\mid\t_A>t)\;,
$$
and starting from $\a$ the hitting time to $A$ is exponential with rate
$1/\E_\a\t_A$.
 If the set $A$ is such that $R/\E_\pi\t_A$ is small, then the distance between
 the stationary and  quasi-stationary measures is small and, as shown in Theorem 3 of  \cite{AB1},
 \be{th3AB}
 \P_\pi\bigl(\t_A>t\bigr)\;\ge\; \Bigl(1-{R\over\E_\a\t_A}\Bigr)\, \exp\Bigl\{-{t\over\E_\a\t_A}\Bigr\}\;.
 \ee

 The proofs of  these   results are based on the property of \emph{complete monotonicity} derived from the spectral decomposition of the transition matrix restricted to the complement of $A$.

 While reversibility has been crucially exploited in all these proofs, there exist  another representation for first hitting times, always due to \cite{B3},
based on the famous interlacing eigenvalues theorem of linear algebra.  This representation has been very recently re-derived, using a different approach, \cite{FL} who use it  to the generalize the results to {\bf non reversible} chains.

While the previous formulas are very revealing, they are restricted to the situations in which the ``good'' reference state is the actual equilibrium measure.

The common feature of these approaches is their central use of the invariant measure both as a reference and to compute the corrections to exponential laws.  As this object is generally unknown and hard to control, the resulting theories lead to hypotheses and criteria not easy to verify.

An higher level of generality is achieved by the \emph{martingale approach} recently applied in \cite{BLM} to obtain results comparable to those in theorem \ref{t1} below. In this reference, however, exponential laws are derived for visits to \emph{rare sets} ---that is, sets with asymptotically small probability.  In metastability or reliability theory this corresponds to visits to the \emph{saddles} mediating between good and bad or between stable and metastable states.
As mentioned above, we recall that the approach proposed in the present paper does not require the determination of these saddle states, and exponential laws
are derived also for visit to sets $G$  including the stable state. Moreover, in this work we concentrate on \emph{recurrence hypotheses} defined purely in terms of the stochastic evolution with no reference to an eventual invariant state.

\subsection{Metastability: A few key settings.}

Having in mind different asymptotic regimes, many different definitions of metastable states have been given in the literature.
These notions, however, are not completely equivalent as they rely on different properties of hitting and escape times.
This state of affairs makes direct comparisons difficult and may lead to confusion regarding applicability of the different theories to new problems.

Since metastability is always associated with a particular asymptotic regime, the results given in the literature are always given in asymptotic form.

In order to understand the reasons behind the different notions of metastability, it is useful to survey the main situations where metastability has been studied.

The simplest case is when the system recurs in a single point.
The main asymptotic regimes that fall into this class are:

\begin{itemize}
\item Finite state space in the limit of vanishing transition probabilities.
Typical examples are lattice systems, with short range interaction, Glauber \cite{AC,BM,CL,CN2012,CO,KO,KO2,NO,NeS,NevSchbehavdrop,S1,S2}(or Kawasaki\cite{BHN,GOS,HNOS,HNT1,HNT2,HNT3,dHOS,NOS}, or parallel\cite{BCLS,CaTr,CN,CNS,CN2012,NS1,Trouve1,Trouverl}) dynamics, in the limit of vanishing temperature.
In this regime, the transition probabilities between neighbour points $x$ and $y$ have the form $P(x,y)=\exp(-\b \D H(x,y))$, where $\beta \to \infty$ is the inverse temperature and $\D H(x,y)$ is the energy barrier between   $x$ and $y$.\\
In this regime, the escape pattern can be understood in terms of the energy landscape (more generally, in terms of the equilibrium measure).
The mean exit time scales as $\exp (-\b \Gamma)$, where $\Gamma$ is the energy barrier between the metastable and the stable configurations.

\item Finite transition probabilities in the limit of diverging number of steps to attain stability.
The typical example of this regime regards mean field models \cite{BEGK01,CGOV}.  Indeed, these systems can be mapped into a low-dimensional Markov process recurring in a single point and where metastability can be described in terms of free--energy landscape.
The mean exit time scales as $\exp (-n \Delta)$,
where $n \to \infty$ is the volume and $\Delta$ is the free-energy barrier divided by the  temperature.

The kinetic Ising model at finite temperature in the limit of vanishing magnetic field $h$ in a box of side-length $1/h$ is another example of this class (see \cite{SS}), since the critical droplet becomes larger and larger as $h$ tends to $0$.

Many toy models, including the one--dimensional
model discussed in section \ref{abc}, pertain to this regime.
\end{itemize}

In the rest of the subsection we will discuss some cases where metastability can be described in terms of entropic corrections of the regimes above.

For instance, suppose to have many independent copies of a system in the regime of finite state space and vanishing transition probabilities, and that the target event is the hitting to a particular set in any of the copies.
Physical phenomena of this kind are bolts (discharge of  supercharged condensers) and "homogeneous nucleation" in short range thermodynamic systems (i.e. the formation of the first critical nucleus in a large volume) \cite{BHS,dHOS,GHNOS09,DS2,SS}
Since the target event can take place in any of the subsystems, the hitting time is shortened with respect to a single subsystem and, in general,  its law changes.
This regime is the physical motivation behind our notion of "early behavior" (see definition \ref{EE} below).

	A wonderful example where the entropic correction is related to fine details of the dynamics is given by the "nucleation and growth models" (see e.g. \cite{DS1,DS2,CeMa,MO1,MO2,SS}), where the transition to stability is driven by the formation, the growth and eventually the coalescence of critical droplets.  In these systems the mean relaxation time is the sum of the "nucleation time" in a "critical volume", which is often exponentially distributed, and the "travel time" needed for two neighbor droplets to grow and coalesce (which, at least in some systems, is believed to have a cut-off behavior).

For general Markov chains, however, metastability cannot be understood in terms of the invariant measure landscape.
The system is trapped in the metastable state both because of the hight of the saddles and the presence of bottlenecks and there is not a general recipe to analyze this situation.
The results given in this paper allow to deal with the cases where the system recurs in a single point.
In a forthcoming paper, we will deal with the general case, where recurrence to the "quasi--stationary measure" (or to a sufficiently close measure) is used.

In section \ref{comphyp}, we compare the different definitions of metastability that have been given in the literature for different asymptotic regimes.

\subsection{Goal of the paper}

The goal of this paper is twofold.  First, we develop an overall approach to the study of exponential hitting times which, we believe, is at the same time general, natural and computationally precise.  It is general because it does not assume reversibility of the process, does not assume that the hitting is to a rare event  and does not assume a particular starting measure. It is natural because it relies on the most universal hypothesis for exponential behavior ---recurrence--- without any further mathematical assumptions like complete monotonicity or other delicate spectral properties of the chain.  Furthermore, the proofs are designed so to follow closely physical intuition.  Rather than resorting to powerful but abstract probabilistic or potential theoretical theorems, each result is obtained by comparing escape and recurrence time scales and decomposing appropriately the relevant probabilities.  Despite its mathematical simplicity, the approach yields explicit bounds on the corrections to exponentiality, comparable to Lemma 7 of \cite{AB1} but without assuming initial equilibrium measure.

A second goal of the paper is to compare and relate metastability conditions proposed in the literature \cite{BovMetaPTA,BEGK01,MNOS}.  Indeed, different authors rely on different definitions of metastable states involving different hypotheses on hitting and escape times.  The situation is particularly delicate for evolutions of infinite-volume systems, whose treatment depends on whether and how relevant parameters (temperature, fields) are adjusted as the thermodynamic limit is taken.  We do a comparative study of these hypotheses to eliminate a potential source of confusion regarding applicability of the different theories to new problems.

A further contribution of our paper is our notion of \emph{early asymptotic exponential behavior} (Definition \ref{exptail}) to control the exponential behavior on a time scale asymptotically smaller than the mean-time scale.  This notion is particularly relevant for systems with unbounded state space where it leads to estimations of probability density functions.  Furthermore, as discussed below, this strong control of exponentiality at small times is important to control infinite-volume systems in which the nucleations leading to exit from metastability can happen anywhere in the volume. We provide natural sufficient conditions on recurrence times for this early exponentiality to hold.

 The main limitation of our approach ---shared with the majority of the metastability literature--- is the assumption that recurrence refers to the visit to a particular  configuration. Actually this particular configuration $x_0$ can be chosen quite arbitrarily in the "metastable well". A more general treatment, involving extended metastable measures for which it is not possible to speak of metastable well, is the subject of current research (\cite{FMNScSo}).

We conclude this section with the outline of the paper.
In section 2 we give main definitions, results on the exponential behavior and compare different hypothesis used in the literature with the ones used in this paper. In section 3 we give some key lemmas that are used in section 4 to prove the main theorems about the exponential behavior.
In section 5 we prove the results about the comparison of different hypothesis and we give an example that depending on the values of parameters fulfills different hypothesis to show that in general they are not all equivalent.  The appendix contains computations of quantities needed  in the example.

\section{Results}
\subsection{Models and notation}

We consider a family of discrete  time irreducible Markov chains with transition matrices
$P^{(n)}$ and invariant measures $(\pi^{(n)}, \, n\ge 1)$ on finite state spaces   $(\cX^{(n)}, \, n\ge 1)$.
A particular interesting case is the infinite volume asymptotics $\lim_{n\to\infty} |\cX^{(n)}|=\infty$.

We use scriptless symbols $\P(\cdot)$ and $\E(\cdot)$ for probabilities and expectations while keeping the parameter $n$ and initial conditions as labels of
events and random variables.  In particular $X^{(n),x}=\bigl(X^{(n),x}_t\bigr)_{t\in\mathbb N}$ denotes the chain
starting at $x\in\cX^{(n)}$, and the hitting time of a set
$F^{(n)}\subset \cX^{(n)}$ is denoted by
\be{hitting}
\t^{(n),x}_{F^{(n)}}=\inf \Bigl\{t\ge 0 :\; X^{(n),x}_t\in F^{(n)}\Bigr\}
\ee

For guidance we
use uppercase latin letters for sequences of diverging positive constants,
and lowercase latin letters for numerical sequences converging to zero when $n\to\infty$.
The small-o notation $ o_n(1)$ indicates a sequence of functions going uniformly to $0$ as $n\to\infty$.
The symbol $A_n\succ B_n$ indicates, $B_n/A_n=o_n(1)$.

The notation $X=\bigl(X_t\bigr)_{t\in\mathbb N}$ is used for  a generic chain on $\cX$.
Quantitave results will be given for such a generic chain, while asymptotical results are given
for the sequence  $X^{(n)}$. The shorter  notation without superindex $(n)$ is also used in
proofs where result are discussed for a single choice of $n$.


\subsection{General results on exponential behavior}

The general setup in the sequel comprises some ingredients.  First, a point $x_0$ thought of as a(meta)stable state.  In reversible chains, such a point corresponds to the bottom of an ``energy well" or, more generally, to a given state in the energy well.  In our treatment, it is irrelevant whether it corresponds to an absolute (stable) or local (metastable) energy minimum.
 The second ingredient is a non empty set $G$ of points marking the exit from the ``well''.
 Depending on the application, this set can be formed by exit points, by saddle points, by the basin of attraction of the stable points or by the target stable points.
 The random time  $\t_{G}^{x_0}$,  i.e.,  the first hitting time to $G$ starting from $x_0$,  corresponds therefore to the \emph{exit time or transition time}
 in the metastable setting. We call  $(x_0, G)$ a  {\it reference pair}.

We characterize the scale of return times (renewal times) by means of two parameters:
\bd{kl1}
Let  $R>0$ and $r\in(0,1)$, we say that  a reference pair $(x_0, G)$ satisfies $Rec(R,r)$ if
\be{ric1}
  \sup_{x\in\cX} \P\Bigl(\t_{\{x_0,G\}}^{x} > R\Bigr) \;\le\; r\;.
 \ee
 \ed
We will refer to $R$ and $r$ as the {\it recurrence time} and {\it recurrence error} respectively.
The hitting time to $\{x_0,G\}$ is one of the key ingredients of our renewal approach.
\bd{bacino}
Given a reference pair  $(x_0, G)$ and $r_0\in (0,1)$, we define the basin of attraction of $x_0$ of confidence level $r_0$
the set
\be{def.bac}
B(x_0,r_0):=\{x\in\cX;\; \P(\t^x_{\{x_0,G\}}=\t^x_{x_0})>1-r_0\}
\ee
\ed

\bt{t0} Consider a reference pair $(x_0, G)$, with $x_0\in\cX$, $G\subset \cX$, such that $Rec(R,r)$ holds with
$R<T:=\E \t_{G}^{x_0}$, with $\e:= {R\over T}$ and $r$ sufficiently small.  Then, there exist functions $C(\epsilon,r)$ and $\lambda(\epsilon,r)$ with
\be{eq:rconv}
C(\epsilon,r)\; ,\;
\lambda(\epsilon,r)\; \longrightarrow\;0\quad\mbox{as}\quad
\epsilon, r\,\to\, 0\;,
\ee
such that
\be{expquant}
\Big|\P\Big({\t_{G}^{x_0}\over T} >t\Big)-e^{-t}\Big|\;\le\;
C\,e^{-(1-\lambda)\,t}
\ee
for any $t>0$.  Furthermore, there exist a function $\widetilde C(\epsilon,r,r_0)$ with
\be{eq:rcc}
\widetilde C(\epsilon,r,r_0)\to 0 \quad\mbox{as}\quad
\epsilon, r, r_0 \,\to\, 0\;,
\ee
such that, for any $z\in B(x_0,r_0)$,
\be{ct0}
\Big|\P\Big({\t_{G}^{z}\over T} >t\Big)-e^{-t}\Big|\;\le\;
\widetilde C\,e^{-(1-\lambda)\,t}
\ee
\et

\begin{rems}
As already discussed in the Introduction, similar results are given in the literature both in the
field of first hitting to rare events and in the field of escape from metastability.
We have to stress here  that we are not assuming reversibility and we have quite general assumptions
on starting condition.
\end{rems}

Metastability studies involve sequences of Markov chains on a sequence of state spaces $\cX^{(n)}$.

\bd{hpG}\label{def2.4}
A sequence of reference pairs $\bigl(x_0^{(n)},G^{(n)}\bigr)$, with $x_0^{(n)}\in\cX^{(n)}$, $\emptyset \neq G^{(n)}\subset \cX^{(n)}$ satisfies Hypothesis Hp.$G(T_n)$ for some increasing positive sequence $(T_n)$ if there exist sequences $r_n=o_n(1)$ and $R_n\prec T_n$ such that
  the recurrence property $Rec(R_n,r_n)$ holds in the following sense:
   \be{ric1.1}
 \sup_{x\in\cX^{(n)}}   \P\Bigl(\t_{\{x_0^{(n)},G^{(n)}\}}^{(n),x} > R_n\Bigr) \;\le\; r_n\;.
  \ee
\ed
\bd{bacinoasym}
Given a sequence of reference pairs  $(x_0^{(n)}, G^{(n)})$ and
a sequence  $r^{(n)}_0\to 0$, we define the basin of attraction of $x^{(n)}_0$ of confidence level $r^{(n)}_0$
the set
\be{def.bacasym}
B(x^{(n)}_0,r^{(n)}_0):=\{x\in\cX^{(n)};\; \P(\t^x_{\{x_0^{(n)},G\}}=\t^x_{x_0^{(n)}})>1-r^{(n)}_0\}
\ee
\ed

We have the following exponential behaviour for the first hitting time to $G^{(n)}$:

 \bt{t1}
 Consider a sequence of reference pairs $\bigl(x_0^{(n)},G^{(n)}\bigr)$ with mean exit times
\be{defTE}
  T_n^E:=\E\Bigl(\t_{G^{(n)}}^{(n),x_0^{(n)}}\Bigr)
\ee and $\zeta$-quantiles information time
\be{defQzeta}
Q_n(\zeta):=\inf\Bigl\{k\ge 1\,:\,  \P\Bigl(\t_{G^{(n)}}^{(n),x_0^{(n)}}  \le k\Bigr)\ge 1-\zeta\Bigr\}.
\ee
Then,
\begin{itemize}
\item[(I)] If Hp.$G(T_n^E)$ holds:
 \begin{itemize}
 \item[(i)]
$\t_{G^{(n)}}^{(n),x_0^{(n)}}/T_n^E$   converges in law  to an $\exp(1)$ random variable,
that is,
  \be{th1.2}
        \lim_{n \to \infty} \P\Bigl(\t_{G^{(n)}}^{(n),x_0^{(n)}}  > t \,T_n^E\Bigr) \;=\; e^{-t}\;.
  \ee
\item[ii)] Furthermore,
  \be{th1.2bac}
        \lim_{n \to \infty}\sup_{x\in B(x^{(n)}_0,r^{(n)}_0)} \Big| \P\Bigl(\t_{G^{(n)}}^{(n),x}  > t \,T_n^E\Bigr) \;-\; e^{-t}\Big|=0
  \ee
\end{itemize}
\item[(II)] If Hp.$G(Q_n(\zeta))$ holds:
\begin{itemize}
\item[(i)]
$\t_{G^{(n)}}^{(n),x_0^{(n)}} /Q_n(\zeta)$   converges in law  to an $\exp(-\ln \zeta)$ random variable, that is,
\be{th1.2.1}
        \lim_{n \to \infty} \P\Bigl(\t_{G^{(n)}}^{(n),x_0^{(n)}}  > t \,Q_n(\zeta)\Bigr) \;=\; \zeta^{t}\;.
  \ee
\item[(ii)] The rates converge,
\be{E/T}
   \lim_{n \to \infty}{Q_n(\zeta)\over T_n^E }=- \ln \zeta
  \ee
   \end{itemize}
   \end{itemize}
  \et

A popular choice for the parameter $\zeta$ in the quantile is $e^{-1}$. In this case all rates are equal to $1$.  For our purposes all choices will work.

Under weaker hypotheses on $T_n$, we can prove exponential behavior in a  weaker form:
\bc{ct1}
If $T_n$ is a sequence of times such that
$$
\exists \, \zeta<1 \,:\, P\Bigl(\t_{G^{(n)}}^{(n),x_0^{(n)}} >  \,T_n\Bigr) \ge  \zeta
\mbox{ uniformly in } n
$$
and if Hp.$G(T_n)$ holds, then the sequence $\t_{G^{(n)}}^{(n),x_0^{(n)}} /\E\Bigl(\t_{G^{(n)}}^{(n),x_0^{(n)}}\Bigr)$
converges in law to an exp(1) random variable. Indeed in this case we have $T_n\le Q_n(\zeta)$ so that Hp.$G(T_n)$
implies Hp.$G(Q_n(\zeta))$.
\ec

Notice that
recurrence to a single point is not necessary to get exponential behavior of the hitting time. A simple example where we get exact exponential behavior independently of the initial distribution is when the one-step transition probability $P_{x,G}$ is constant for $x\in G^c$.

\medskip

As explained in the introduction there are situations, e.g., when treating metastability in large volumes, that call for more detailed information on a short time scale on hitting times.

\bd{EE}
A random variable $\theta$ has early exponential behaviour at scale $S\le\E \theta$ with rate $\a$, if for each $k$ such that
$kS\le\E\theta$ we have
\be{eequant}
\Big|{\P\bigl(\theta\in(  k S,(k+1)S]\bigr)\over  \P(\theta >   S)^k\,\,\P(\theta \le  S)}-1\Big|<\a.
\ee
We denote this behavior by $EE(S,\a)$.
\ed

\br{ree}
A remark on the difference between the notion of early exponential behaviour and
the exponential behaviour given by Theorem \ref{t0} is necessary.
Early exponential behaviour controls the distribution only on short times, $kS<\E\theta$, while equation (\ref{expquant}) holds for any $t$. However on the first part of the distribution EE$(S,\a)$ can give a more detailed control on the density of the distribution. More precisely if $\t_G^{x_0}$ is EE$(S,\a)$ with $\a$ small, we can obtain estimates
on the density $f(t)$ of
${\theta\over\E\theta}$,  equivalent to the results obtained in \cite{AB1}, Lemma 13, (a),(b).
Indeed
$$
\P\bigl(\t_G^{x_0}\in(  k S,(k+1)S]\bigr)\;=\;e^{-\l k}(1-e^{-\l})(1+a_k)
$$
where $\lambda := -\ln \P (\t_G^{x_0} >S)$ and
$$
a_k\;:=\;{\P\bigl(\t_G^{x_0}\in(  k S,(k+1)S]\bigr)\over  \P(\t_G^{x_0}>   S)^k\,\,\P(\t_G^{x_0} \le  S)}-1
$$
The absolute value of $a_k$ is bounded  by $\alpha$ if EE$(S,\a)$ holds
 uniformly in $k<\frac{\E\t_G^{x_0}}{S}$. In the case $S\ll \E\t_G^{x_0}$, Lemma \ref{kl3} below implies that
 $\lambda\sim{S\over \E\t_G^{x_0}}$.  Thus,

\begin{eqnarray}
f\Bigl(k{S\over\E\t_G^{x_0}}\Bigr) \,{S\over \E\t_G^{x_0}}&\sim&
\P\Bigl({\t_G^{x_0}\over \E\t_G^{x_0}}\in(  k {S\over \E\t_G^{x_0}},(k+1){S\over\E\t_G^{x_0}}]\Bigr)\nonumber\\
&=&e^{-\l k}(1-e^{-\l})(1+a)\nonumber\\
&\sim& e^{-{S\over \E\t_G^{x_0}} k}\,\Bigl[{S\over \E\t_G^{x_0}}+ o\Bigl({S\over \E\t_G^{x_0}}\Bigr)\Bigr]\,\bigl(1+a\bigr)
\label{eq:new1}
\end{eqnarray}
\er


\bd{exptail}
A family of random variables $(\theta_n)_n$ with $\E\theta_n\to\infty$ for $n\to\infty$,  has an  asymptotic early exponential behaviour at scale $(S_n)_n$ if for every integer $k$

\be{expbeh}
\lim_{n\to\infty}{\P\bigl(\theta_n\in(  k S_n,(k+1)S_n]\bigr)\over  \P(\theta_n >   S_n)^k\,\,\P(\theta_n \le  S_n)}=1
\ee

\ed

\br{ree2}
The notion of $EE(S_n,\a)$ is interesting for $\a$ small and when $S_n$ is
asymptotically smaller than $\E(\theta_n)$ so that $\P(\theta_n\leq   S_n) \to 0$ as
$n\to \infty$.  The sharpness condition  \eqref{expbeh}
controls the smallness of these probabilities.  In particular it implies that
\be{exptail2}
\frac{\P(\theta_n\leq  k S_n)}{1-\P(\theta_n >   S_n)^k} \;\tendn\; 1
\ee
\er

The following theorem determines convenient sufficient conditions for early exponential behaviour at scale.

\bt{thEEquant}
Given a reference pair $\bigl(x_0,G\bigr)$,  satisfying   $Rec(R,r)$ with
$0<R<T:=\E \t_{G}^{x_0}$.
Define $\e:= {R\over T}$ and suppose $\e$ and $r$ sufficiently small. Then $\t_{G}^{x_0}$ has $EE(\eta T,\a)$ with $\alpha=O(\epsilon/\eta)+O(r/\eta)$, for $\eta$ such that $\eta\in (0,1)$, $\epsilon/\eta$ and $r/\eta$ are small enough.  For instance the property holds for
$\eta=\bigl[\max\{\epsilon,r\}\bigr]^\gamma$ with $\gamma<1$ and $\epsilon, r$ small enough.
\et
The following theorem is an immediate consequence of the previous one.
\bt{thEEasym}
Consider a sequence $\bigl(x_0^{(n)},G^{(n)}, T_n\bigr)$, with $x_0^{(n)}\in\cX^{(n)}$, $G^{(n)}\subset \cX^{(n)}$ and $T_n>0$ satisfying Hypothesis Hp.$G(T_n)$, with $r_n\to 0$. Then
 the family of random variables ${\t_{G^{(n)}}^{(n),x_0^{(n)}}}$
 has  asymptotic exponential behavior at every scale $S_n$  such that
 \be{eq:rm100}
 r_n\prec\; \frac{S_n}{T_n}\qquad\hbox{ and } \qquad\frac{R_n}{T_n}\;\prec\; \frac{S_n}{T_n}\;\le\; 1\;.
 \ee
\et


\subsection{\label{comphyp}Comparison of hypotheses}

Many different hypotheses have been used in the literature to prove exponential behavior.
In this section we analyze some of these hypotheses in order to clarify the relations between them.
The notation \ref{def2.4} can be used to discuss different issues, in particular:
\begin{itemize}
\item the hitting problem (where $x_0$ is the maximum of the equilibrium measure and $G$ is a rare set)
\item the exit problem and applications to metastability(where $x_0$ is a local maximum of the equilibrium measure, e.g., the metastable state, and $G$ is either the ''saddle", the basin of attraction of the stable state or the stable state itself)
\end{itemize}

A key quantity, especially in the "potential theoretic approach" is the following:
\bd{LocalTimeMS}
$A\subset \cX$ and $z,x\in\cX$ the local time spent in $x$ before reaching $A$ starting from $z$ is
\be{loctime}
\xi^{z}_{A}(x):=\Bigl|\Bigl\{t<\t^{z}_{A}:\; X^{z}_t=x\Bigr\}\Bigr|.
\ee
\ed

The following hypotheses are instances of Definition \ref{hpG}.
\medskip

\noindent{\bf Hypotheses I}
\begin{eqnarray}\label{HypG}
   &Hp.G^E  \equiv Hp.G(T_n^E)\quad &\hbox{  with  } \quad T_n^E=\E\Bigl(\t_{G^{(n)}}^{(n),x_0^{(n)}}\Bigr) \\
   &Hp.G^{\zeta}   \equiv Hp.G(T_n^{Q^\zeta})\quad &\hbox{  with  } \quad T_n^{Q^\zeta}=\inf\Bigl\{t:\;\P\Bigl(\t^{(n),x_0^{(n)}}_{G^{(n)}}>t\Bigr)\leq\zeta\Bigr\}\\
   &Hp.G^{LT}  \equiv Hp.G(T_n^{LT})\quad &\hbox{  with  } \quad T_n^{LT}=\E\Bigl(\xi^{(n),x_0^{(n)}}_{G^{(n)}}(x_0^{(n)})\Bigr)
\end{eqnarray}

We show in Theorem \ref{t3} below that the first two hypotheses are equivalent for every $\zeta<1$. In particular this shows the insensitivity of metastability studies to the choice of $\zeta$.  For  $\zeta=e^{-1}$, hypothesis $Hp.G^{1/e}$ is equivalent to the ones considered in previous papers (see for instance \cite{MNOS} hypothesis of Theorem 4.15) to determine the distribution of the escape times for general Metropolis Markov chains in finite volume.
 The last hypotheses $Hp.G^{LT}$ is new and it is useful to compare the first two hypothesis with the two hypotheses below.

The next set of hypotheses refer to the following quantity.
\bd{firshittimePT}
Let
\be{tau+}
\widetilde\t^{{(n)},x}_A:=\min\Bigl\{t>0: X^{(n),x}_t\in A\Bigr\}\ee
be the first positive hitting time to $A$ starting at $x$.
\ed

\noindent
Given reference pairs $\{x_0^{(n)}, G^{(n)}\}$, define
\be{roA}
\rho_A(n):=\sup_{z\in \cX^{(n)}\backslash \{x_0^{(n)}, G^{(n)}\}}
{\P\Bigl(\widetilde\t^{{(n)}, x_0^{(n)}}_{G^{(n)}}<\widetilde\t^{{(n)},x_0^{(n)}}_{x_0^{(n)}}\Bigr)
\over \P\Bigl(\widetilde\t^{{(n)},z}_{\{x_0^{(n)}, G^{(n)}\}}<\widetilde\t^{{(n)},z}_z\Bigr)}
\ee
and
\be{roB}
\rho_B(n):=\sup_{z\in \cX^{(n)}\backslash \{x_0^{(n)}, G^{(n)}\}}{\E \t^{{(n)},z}_{\{x_0^{(n)}, G^{(n)}\}}\over \E \t^{(n), x_0^{(n)}}_{G^{(n)}}}\;.\qquad\qquad
\ee
\medskip

\noindent{\bf Hypotheses II}
\begin{eqnarray}\label{HypAB}
   &Hp. A: \quad &\lim_{n\to\infty}|\cX^{(n)}|\rho_A(n)=0\qquad\qquad\qquad\\
   &{Hp.B}: \quad &\lim_{n\to\infty}\rho_B(n)=0\qquad\qquad\qquad
\end{eqnarray}

When $x_0$ is the metastable configuration and $G$ is the stable configuration (more precisely when $\E (\xi^{x_0}_G)<\E (\xi_{x_0}^G)$ )  $Hp. A$ is similar to the hypotheses considered in \cite{BEGK00} while $Hp. B$ is similar to those assumed in \cite{BovMetaPTA}. Theorem 1.3 in \cite{BEGK02} shows that, under hypothesis $Hp. A$ and reversibility, $\t^{x_0}_G / \E (\t^{x_0}_G)$ converges to a mean $1$ exponential variable.


Our last theorem establishes the relation between the previous six hypotheses.
\bt{t3}
The following implications hold:
\be{eq:implic}
Hp.A\;\Longrightarrow \;Hp.G^{LT}\;\Longrightarrow\; Hp.G^{E}\;\Longleftrightarrow \;  Hp.G^{\zeta}\;  \Longleftrightarrow \;Hp.B
\ee
for any $\zeta<1$.  Furthermore, the missing implications are false.
\et

\begin{figure}
\centering
\includegraphics[width=0.5\linewidth]{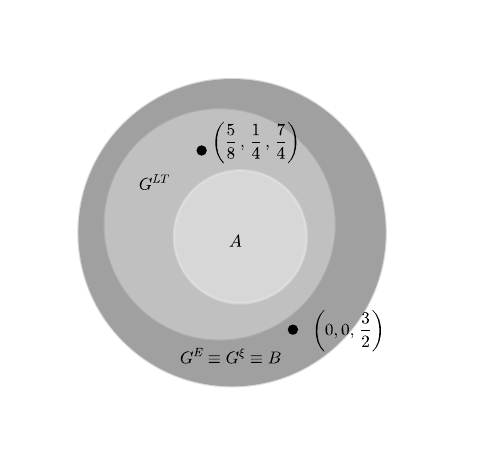}
\caption{Venn diagram for the conditions in theorem \ref{t3}. The two points correspond to particular choices of the parameters in the "abc model" stated in section \ref{abc}. INGRANDIRE}
\label{fig:Venn2}
\end{figure}

These relations are summarized in Figure \ref{fig:Venn2}.

\begin{rems}\
\bi
\item[(a)] Theorem \ref{t1} holds with either hypothesis  $Hp.B$ of the originally stated $Hp.G^E$ or $Hp.G^{\zeta}$.
\item[(b)] The first equivalence shows, in particular, that all hypotheses $Hp.G^{\zeta}$ with $\zeta<1$ are mutually equivalent.
\item[(c)] Typically, in metastable systems, for a given target set $G$ there are many possible choices for the point $x_0$ to form a reference pair that verifies Hp. B.\\
Let
$\mathcal{M}_{G}^{\epsilon}:=\left\{ x\ ;\;\sup_{z\not\in\left\{ x,G\right\} }\mathbb{E}(\tau_{\left\{ x,G\right\} }^{z})/\mathbb{E}(\tau_{G}^{x})<\epsilon\right\} $
be the set of all points that toghether with $G$ form a reference pair
according to Hp.B.\\
If $x,\; y\in\mathcal{M}_{G}^{\epsilon}$, then
\be{equivT}
\frac{1}{1+\epsilon}\le\frac{\mathbb{E}(\tau_{G}^{x})}{\mathbb{E}(\tau_{G}^{y})}\le1+\epsilon
\ee
and
\be{Pvisit}
\mathbb{P}\left(\tau_{x}^{y}<\tau_{G}^{y}\right)\ge1-2\epsilon\ \ \ ;\ \ \ \mathbb{P}\left(\tau_{y}^{x}<\tau_{G}^{x}\right)\ge1-2\epsilon
\ee
Indeed, for any $x,y,z\in\mathcal{X}$,
\begin{eqnarray}
\mathbb{E}(\tau_{z}^{y})&=&\mathbb{E}(\tau_{\left\{ x,z\right\} }^{y})+\mathbb{E}(\tau_{z}^{y}-\tau_{\left\{ x,z\right\} }^{y})\nonumber\\
&=&\mathbb{E}(\tau_{\left\{ x,z\right\} }^{y})+\mathbb{E}\left((\tau_{z}^{y}-\tau^y_x){1}_{\tau_{x}^{y}<\tau_{z}^{y}} \right) \nonumber\\
&=&\mathbb{E}(\tau_{\left\{ x,z\right\} }^{y})+\mathbb{E}(\tau_{z}^{x})\mathbb{P}\left(\tau_{x}^{y}<\tau_{z}^{y}\right),
\label{eq:renewT}
\end{eqnarray}
where we used strong Markov property at time $\tau_{x}^{y}$ in the last equality.\\
When $x\in\mathcal{M}_{G}^{\epsilon}$,
\be{above}
\mathbb{E}(\tau_{G}^{y})\le\mathbb{E}(\tau_{G}^{x})\left(\epsilon+1\right)
\ee
Hence, if both $x$ and $y$ are in  $\mathcal{M}_{G}^{\epsilon}$,
we get  the (\ref{equivT}). From (\ref{eq:renewT}), by using (\ref{equivT}) and $x\in\mathcal{M}_{G}^{\epsilon}$, we get
\[
\mathbb{P}\left(\tau_{x}^{y}<\tau_{G}^{y}\right)=\frac{\mathbb{E}(\tau_{G}^{y})-\mathbb{E}(\tau_{\left\{ x,G\right\} }^{y})}{\mathbb{E}(\tau_{G}^{x})}\ge\frac{1}{1+\epsilon}-\epsilon,
\]
by  symmetry,  (\ref{Pvisit}) follows.
\item[(d)] A well-known case is that of finite state space under Metropolis dynamics in the limit of vanishing temperature.
Lemma 3.3 in \cite{BM} states that if $G$ is the absolute minimum of the energy function and $x_0$ is the deepest local minimum, then the reference pair $(x_0,G)$ verifies Hp. A.
\ei
\end{rems}



\section{Key Lemmas}

The proofs of the theorems presented in this paper are based on some quite simple results
on the distribution of the random variable $\t^{x_0}_G$.  In fact, the central argument is that condition $Rec(R,r)$ implies that
renewals ---that is, visits to $x_0$--- happen at a much shorter time scale than visits to $G$.  This is expressed through the behavior of the following random times related to recurrence:
for  each deterministic time $u>0$ let
\be{eq:r3}
\t^*(u)\;:=\;\inf\Bigl\{{s \ge u \tc X_s\in\{{x_0,G}}\}\Bigr\}
\ee
If $Y\sim Exp(1)$ then the following factorization obviously holds
$$
\P(Y>t+s)\;=\;\P(Y>t)\,\P(Y>s).
$$
Next  lemma controls --for $z=x_0$-- this factorization property on a generic time scale $S>R$.
\bl{kl2}
 If $(x_0,G)$ satisfies $Rec(R,r)$, then for any $z\in \cX$,
   $S>R$, $t>0$ and $s>{R\over S}$
\be{eq:r2}
\P\Bigl({\t^{z}_G > (t+s) S} \Bigr)
\begin{array}{l}
\ge \;
\Bigl[ \P\Bigl({\t^{z}_G > t S + R}\Bigr)-r\,\P\Bigl({\t^{z}_G > t  S}\Bigr) \Bigr]\,
  \P\Bigl({\t^{x_0}_G > s  S}\Bigr)\\
\le\Bigl [{
   \P\Bigl({\t^{z}_G > t  S - R}\Bigr)+
   r}\Bigr]\P\Bigl({\t^{x_0}_G > s  S}\Bigr)\,
 \;.
   \end{array}
\ee
\el

\bpr
We start by decomposing according to the time $ \t^*(t S)$ to get
\begin{eqnarray}\label{t*1}
 \lefteqn{
  \P\Bigl({\t^{z}_G > (t+s) S} \Bigr) } \nonumber\\[5pt]
  &=&
  \P\Bigl({\t^{z}_G > (t+s) S\; ;\; \t^*(t S) \le t S + R}\Bigr)\nonumber\\
  &&\qquad+\;
  \P\Bigl({\t^{z}_G > (t+s) S \; ;\; \t^*(t S) >t S + R}\Bigr)\\[5pt]
 &=& \sum_{u=0}^{R}
  \P\Bigl({\t^{z}_G > (t+s)  S\; ;\; \t^*(t S) = t  S + u}\Bigr)\nonumber\\
   &&\qquad +\;
  \sum_{x \in  \{x_0,G\}^c }
  \P\Bigl({\t^{z}_G > (t+s)  S \; ;\; X_{t  S}=x \; ;\; \t^*(t S) >t  S + R}\Bigr)
  \nonumber
\end{eqnarray}
We now use Markov property at time $t S+u$ in the first sum, together with the fact that  $\t^{z}_G>\t^*(tS)$ implies $X_{\t^*(tS)}=x_0$.
In the second sum we use Markov property at instant $t  S$.  This yields
\begin{eqnarray}\label{t*1tre}
  \lefteqn{\P\Bigl({\t^{z}_G > (t+s)  S}\Bigr) }\nonumber\\
  &=&\sum_{u=0}^{R}
  \P\Bigl({
      \t^*(t S) = t  S + u \; ;\; \t^{z}_G > t  S + u}\Bigr)\,
  \P\Bigl({\t^{x_0}_G > s  S - u}\Bigr)   \\
   &&\qquad +\;
  \sum_{x \in  \{x_0,G\}^c }
  \P\Bigl({\t^{z}_G > t  S        \; ;\; X_{t  S}=x}\Bigr)\,
  \P\Bigl({\t^x_G > s  S \; ;\; \t^x_{\{x_0,G\}} >  R }\Bigr)\;.\nonumber
\end{eqnarray}
\smallskip
This identity will be combined with the elementary monotonicity bound with respect to  inclusion
\be{monot}
\{\t^{x_0}_G>t_1\}\supseteq \{\t^{x_0}_G>t_2\}\quad\hbox{ for } t_1\le t_2\;,
\ee

To get the lower bound we  disregard the second sum in the right-hand side of \eqref{t*1tre} and bound the first line through the monotonicity bound \eqref{monot} for $t_1=tS+u$ and $t_2=tS+R$.  By condition \eqref{ric1} we obtain:
\begin{eqnarray}\label{t*3}
\lefteqn{ \P\Bigl({\t^{z}_G > (t+s)  S}\Bigr )\ge
             \P\Bigl({\t^{z}_G > t  S + R\; ;\;
              \t^*(tS) \le t  S+ R}\Bigr)\,
  \P\Bigl({\t^{x_0}_G > s  S}\Bigr) }  \nonumber\\[8pt]
  &=&
 \Bigl[{\P\Bigl({\t^{z}_G > t  S+ R}\Bigr)-
          \P\Bigl({\t^{z}_G > t  S + R\; ;\;
              \t^*(tS) > t  S+ R}\Bigr)}\Bigr]\,
  \P\Bigl({\t^{x_0}_G > s  S}\Bigr)
\\[8pt]
  &\ge&  \Bigl[{
  \P\Bigl({\t^{z}_G > t  S + R}\Bigr)-  r \,\P\Bigl({\t^{z}_G > t  S}\Bigr)}\Bigr]\,
  \P\Bigl({\t^{x_0}_G > s  S}\Bigr)\nonumber
\end{eqnarray}

To get the upper bound in \eqref{eq:r2} exchange $s$ with $t$ in
(\ref{t*1tre}) and  use again the monotonicity \eqref{monot} to bound  $\P({ \t^{x_0}_G > t  S- u})\leq \P({ \t^{x_0}_G > t  S - R})$
in the first sum in the right-hand side.
For second sum we use condition \eqref{ric1}.  This yields
\begin{eqnarray}\label{t*1due}
\lefteqn{ \P\Bigl({\t^{z}_G > (t+s)  S}\Bigr)}
 \nonumber\\[8pt]
  &\le&  \P\Bigl({
      \t^*(s S)\leq s S+ R  \; ;\; \t^{z}_G > s  S
}\Bigr)\,
  \P\Bigl({\t^{x_0}_G > t  S - R}\Bigr)+
  \P\Bigl({\t^{z}_G > s  S} \Bigr)  \,    r
 \\[8pt]
&\le&  \P\Bigl({\t^{z}_G > s  S}\Bigr)
 \Bigl [{
   \P\Bigl({\t^{x_0}_G > t  S - R}\Bigr)+
   r}\Bigr]\;.\nonumber
\end{eqnarray}
\medskip
\epr

From this ``almost factorization'',  we can control the distribution
of $\t^{x_0}_G $ on time scale $R$. Indeed the following two
results are easy consequences of this factorization.

\bc{cor-fatt}
If $(x_0,G)$ satisfies  $Rec(R,r)$
and $S$ is such that  $R<S<T:=\E\t^{x_0}_G$, then,
\be{iter<}
\P(\t^{x_0}_G>Sk)\le \Big[\P(\t^{x_0}_G>S-R) +r\Big]^k
\ee
\be{iter>}
\P(\t^{x_0}_G>Sk)\ge  \Big[\P(\t^{x_0}_G>S+R) -r\Big]^k
\ee
 for any integer  $ k >1$.
\ec
This Corollary  immediately follows by an iterative application of Lemma \ref{kl2}
with $t=1$ and $s=k-1$.

\bl{kl3}
If $(x_0,G)$ satisfies $Rec(R,r)$,  and $S$ is such that  $R<S<T:=\E\t^{x_0}_G$,
then
\be{Pt<R}
\P\bigl(\t^{x_0}_G\le S\bigr)\;\le\; {S+R\over T}+r
\ee
and
\be{Pt>R}
 \frac{1}{1+\frac{T}{S-2R}}-r\;\le\; \P\bigl(\t^{x_0}_G\le S\bigr)\;.
 \ee
As a consequence,
\be{eq:r0}
\Big|\P\bigl(\t^{x_0}_G\le S\bigr)-{S\over T}\Big| \;<\; 2{R\over T} +\Big({S\over T}\Big)^2+r\;.
\ee
\el

\bpr
Let us denote $m$ the integer part of $S/R$, that is the integer number such that
\be{eq:rr2}
S-R\;<\; mR\;\le\; S\;.
\ee

\noindent
\emph{Proof of the upper bound.} We bound the mean time in the form

\begin{eqnarray}
\label{mean}
\E\t^{x_0}_G=\sum_{t=0}^\infty\P(\t^{x_0}_G>t)
&=&\sum_{k=0}^\infty\sum_{i=(m+1)Rk}^{(m+1)R(k+1)} \P\bigl(\t^{x_0}_G>i\bigr)\\[8pt]
 &\le& (m+1) R \sum_{k=0}^\infty \P\bigl(\t^{x_0}_G>(m+1)Rk\bigr)\;.\nonumber
\end{eqnarray}
The last line is due to the monotonicity property (\ref{monot}).  Therefore, by \eqref{iter<},
\begin{eqnarray}
\label{mean.1}
\E\t^{x_0}_G&\le& (m+1) R \sum_{k=0}^\infty \Big[\P\bigl(\t^{x_0}_G>m \, R\bigr) +r\Big]^k
\nonumber\\[8pt]
&=& (m+1) R \sum_{k=0}^\infty  \Big[1-\P\bigl(\t^{x_0}_G\le mR\bigr) +r\Big]^k\;.
\end{eqnarray}
If $\P\bigl(\t^{x_0}_G\le m R\bigr)\le r$ the bound (\ref{Pt<R}) is trivially satisfied.  Otherwise
the power series converges yielding
$$
\E\t^{x_0}_G
\; \le\; \frac{(m+1) R}{\P\bigl(\t^{x_0}_G\le mR\bigr) -r}
$$
and, thus, by \eqref{eq:rr2},
\be{me+r}
 \P\bigl(\t^{x_0}_G\le S\bigr)\le\P\bigl(\t^{x_0}_G\le mR\bigr)
 \le {(m+1)R\over \E\t^{x_0}_G}+r
 \le \frac{S+R}{T}+r\;.
\ee

\noindent
\emph{Proof of the lower bound.} The argument is very similar, but resorting to the bound
\begin{eqnarray}
\label{mean.10}
\E\t^{x_0}_G&=&\sum_{k=0}^\infty\sum_{i=(m-1)Rk}^{(m-1)R(k+1)} \P\bigl(\t^{x_0}_G>i\bigr)
\nonumber
 \ge (m-1) R \sum_{k=0}^\infty \P\bigl(\t^{x_0}_G> (m-1)R(k+1)\bigr)\nonumber\\[8pt]
 &\ge& (m-1) R \sum_{k=0}^\infty \P\bigl(\t^{x_0}_G> (S-R)(k+1)\bigr)
 \ge (m-1)R \sum_{k=1}^\infty  \Big[\P\bigl(\t^{x_0}_G>S\bigr) -r\Big]^k
 \;.\nonumber
\end{eqnarray}
The second and third inequalities follow from monotonicity, while the last one is due to the lower bound \eqref{iter>}.
The power series is converging because $\bigl|\P\bigl(\t^{x_0}_G>S\bigr) -r\bigr|<1$ since $r\in(0,1)$.  Its sums yields
$$
\frac{\E\t^{x_0}_G}{(m-1)R}
\; \ge\;  \frac{1}{\P\bigl(\t^{x_0}_G\le S\bigr)+r}-1\;,
$$
so that
\be{eq:rr4}
\P\bigl(\t^{x_0}_G\le S\bigr)\;\ge\; \frac{1}{1+\frac{T}{S-2R}}-r
\ee
in agreement with \eqref{Pt>R}s.
\epr


\br{rkl3}
The error term $r$ becomes exponentially small as the recurrence parameter $R$ is increased linearly.  More precisely, if $Rec(R,r)$ holds then
\be{eq:rr5}
\sup_{x\in\cX}\P\Bigl(\t^x_{\{x_0,G\}}>NR\Bigm| \t^x_{\{x_0,G\}}>(N-1)R\Bigr)
\;\le\; r\;,
\ee
which implies,
\be{recR+}
\sup_{x\in\cX}\P\bigl(\t^x_{\{x_0,G\}}>NR\bigl)\;\le\; r^N\;.
\ee
In particular, for $T\gg R$ we can replace $R$ by $R^+:=NR$.
For this reason we can assume in what follows $r<\e$.

\er

The bounds given in the preceding lemmas, however,  are not enough for our purposes. To control large values of $S$ with respect to $T$ (tail of the distribution), we need to pass from additive to multiplicative errors, that is from bounds on  $\bigl|{\P\bigl(\t^{x_0}_G > S + R\bigr)}-{\P\bigl({\t^{x_0}_G > S}\bigr)}\bigr|$ to bounds on $\P\bigl(\t^{x_0}_G > S + R\bigr)\bigm/\P\bigl({\t^{x_0}_G > S}\bigr)$.  Our bounds will be in terms of the following parameters.

\bd{dcbarc}
Let $c, \bar c$ be
\be{eq:defc}
c\;:=\;\P({\t^{x_0}_G \leq 2 R }) +r
\ee
and
\be{eq:defbarc}
\bar c\;:=
\begin{cases}
\frac{1}{2}-\sqrt{\frac{1}{4}- c}& \text{ if} \quad c\le \frac{1}{4}\\
1  & \text{if} \quad c>\frac{1}{4}
\end{cases}
\ee
\ed
We shall work in regimes where these parameters are small.  Note that, by Lemma \ref{kl3}, hypothesis $Rec(R,r)$ implies that
\be{eq:rr7}
c\;<\;3\,\frac{R}{T}+2 r\;.
\ee
Moreover, in the case $c<1/4$ definition (\ref{eq:defbarc}) is equivalent to
\be{cbarc}
c\;=\;\bar c (1-\bar c) \;<\; \bar c\;.
\ee
%
Our bounds rely on the following lemma, which yields control on the tail density of $\t^{x_0}_G $.

\bl{kl4} Let $(x_0,G)$ be a reference pair satisfying $Rec(R,r)$ and $S>R$,   then, for any $z\in \cX$,
\be{stimaleq}
\P\bigl({\t^{z}_G \le S + R}\bigr)\;\le\;
\;\P\bigl({\t^{z}_G \le S }\bigr)\,\Bigl[1+{c\over \P\bigl({\t^{z}_G \le S }\bigr)}\Bigr]
\ee
Furthermore, if $z\in B(x_0, \bar c -c)$,
\be{eq:stimacebarc}
\P\bigl({\t^{z}_G > S + R}\bigr)\;\ge\;
\P\bigl({\t^{z}_G > S }\bigr)\,[1-c-\bar c]\;.
\ee
\el

\bpr

\noindent
\emph{Proof of (\ref{stimaleq}).}
We decompose
\begin{eqnarray}
\label{eq:rr10}
\P\Bigl(\t^{z}_G \le S + R\Bigr)
&=& \P\Bigl({\t^{z}_G \le S }\Bigr)
+\P\Bigl({\t^{z}_G \in ( S, S+ R) }\Bigr)\nonumber\\[5pt]
&=&\P\Bigl({\t^{z}_G \le S}\Bigr)
+\P\Bigl({\t^{z}_G \in ( S, S+ R), \tau^{*}(S-R)\leq S }\Bigr)\nonumber\\
&&\qquad {}+\P\Bigl({\t^{z}_G \in ( S, S+ R), \tau^{*}(S-R)> S }\Bigr)\;.
\end{eqnarray}
The event $\bigl\{\t^{z}_G \in ( S, S+ R), \tau^{*}(S-R)<S \bigr\}$ corresponds to having a visit to $x_0$ in the interval $[S-R,S]$ followed by a first visit to $G$ is in the interval $(S,S+R)$.  By Markovianness,
\begin{eqnarray}
\label{eq:rr12.0}
\lefteqn{\P\Bigl({\t^{z}_G \in ( S, S+ R), \tau^{*}(S-R)\le S }\Bigr)}\nonumber\\[5pt]
&=&\sum_{\scriptstyle u\in[S-R,S]\atop \scriptstyle v\in(S,S+R]}
\P\Bigl(X_u=x_0\,;\, \t^{z}_G\ge u\,;\, X_i\not\in G, u<i<v\,;\, X_v\in G\Bigr)\\[8pt]
&=& \sum_{\scriptstyle u\in[S-R,S]\atop \scriptstyle v\in(S,S+R]}
\P\bigl(X_u=x_0\,,\, \t^{z}_G\ge u\bigr)\, \P\bigl(\t^{x_0}_G=v-u\bigr)\;.\nonumber
\end{eqnarray}
Hence, by monotonicity,
\be{eq:rr12}
\P\Bigl({\t^{z}_G \in ( S, S+ R), \tau^{*}(S-R)\leq S }\Bigr)
\;\le\; \P\Bigl({\t^{z}_G >S- R }\Bigr)\,\P\Bigl({\t^{x_0}_G \le2R}\Bigr)\;.
\ee
Analogously, Markovianness and monotonicity yield
\begin{eqnarray}
\label{eq:rr13.0}
\lefteqn{\P\Bigl({\t^{z}_G \in ( S, S+ R), \tau^{*}(S-R)>S }\Bigr) }\nonumber\\[5pt]
&=&\sum_{x\not\in \{x_0,G\}}
\P\Bigl(\t^{z}_G\in (S,S+R)\,;\, X_S=x\,;\,X_i\not\in \{x_0,G\}, S-R<i\le S\Bigr)\\[8pt]
&\le& \P\bigl(\t^{z}_G>S-R\bigr)\,\sup_{x\not\in\{x_0,G\}}
\P\bigl(\t^{x}_{\{x_0,G\}}>R\bigr)\;.\nonumber
\end{eqnarray}
Hence, by condition $Rec(R,r)$ and monotonicity,
\be{eq:rr13}
P\Bigl({\t^{z}_G \in ( S, S+ R), \tau^{*}(S-R)> S }\Bigr)\;<\;
\P\bigl(\t^{z}_G>S-R\bigr)\,r\;.
\ee
Inserting \eqref{eq:rr12} and \eqref{eq:rr13} in \eqref{eq:rr10} we obtain
\be{eq:rr11}
\P\Bigl(\t^{z}_G \le S + R\Bigr)\;\le\;\P\Bigl({\t^{z}_G \le S}\Bigr)
+ \Bigl[\P\Bigl({\t^{x_0}_G \le2R}\Bigr) +r\Bigr] \P\bigl(\t^{z}_G>S-R\bigr)\;.
\ee
The bound \eqref{stimaleq} is obtained by recalling \eqref{eq:defc} and by neglecting the last term.
\medskip

\noindent
\emph{Proof of (\ref{eq:stimacebarc}).} If $c>1/4$ there is nothing to prove.
For $c\le 1/4$ we perform, for each $z$, a decomposition similar to \eqref{eq:rr10}:
\begin{eqnarray}
\label{eq:rr14}
\P\Bigl(\t^{z}_G > S + R\Bigr)&=& \P\Bigl({\t^{z}_G > S }\Bigr)
-\P\Bigl({\t^{z}_G \in ( S, S+ R) }\Bigr)\nonumber\\[5pt]
&=&\P\Bigl({\t^{z}_G > S}\Bigr)
-\P\Bigl({\t^{z}_G \in ( S, S+ R), \tau^{*}(S-R)<S }\Bigr)\\
&&\qquad {}-\P\Bigl({\t^{z}_G \in ( S, S+ R), \tau^{*}(S-R)\ge S }\Bigr)\;.
\;\nonumber
\end{eqnarray}
The bounds \eqref{eq:rr12} and \eqref{eq:rr13} yield
\be{eq:rr15.0}
\P\bigl(\t^{x}_G > S + R\bigr) \;\geq\;\P\bigl({\t^{z}_G > S}\bigr)
- \P\bigl({\t^{z}_G >S- R }\bigr)\Bigl[ \P\bigl({\t^{x_0}_G \leq 2R}\bigr)+r\Bigr]
\;\nonumber\\
\ee
which can be written as
\be{equation1of3.4}
\frac{\P\bigl(\t^{z}_G > S + R\bigr)}{\P\bigl({\t^{z}_G > S}\bigr)}\;\geq\; 1 - c\,
\frac{\P\bigl(\t^{z}_G > S - R\bigr)}{\P\bigl({\t^{z}_G > S}\bigr)}\;.
\ee
\medskip

Let us first consider the case in which $S=iR$ for an integer $i$.  Denote
\be{definitionofyi}
y_i\;=\;\frac{\P\bigl(\t^{z}_G > (i+1)R\bigr)}{\P\bigl({\t^{z}_G > i R }\bigr)}\;,
\ee
so condition (\ref{equation1of3.4}) becomes
\be{equation2of3.4}
y_i\;\geq\; 1 -
\frac{c}{y_{i-1}}
\ee
The proposed inequality (\ref{eq:stimacebarc}) follows from the following

\noindent
\empbf{Claim:}
\be{equation3of3.4}
y_{i}\;>\; 1 -\bar c
\ee
 We prove this by induction.  For $i=0$, we first notice that
\begin{eqnarray*}
\P\bigl(\t^{z}_G \le 2 R\bigr) =
 \P\bigl(\t^{z}_G \le 2 R\,,\, \t^{z}_{x_0}\ge R\bigr)
 +\P\bigl(\t^{z}_G \le 2 R\,,\, \t^{z}_{x_0}<R\bigr)
 \le \P\bigl(\t^{z}_{\{x_0,G\}}\neq \t^z_{x_0}\bigr)\\
 + \P\bigl(\t^{z}_{\{x_0,G\}}= \t^z_{x_0} \geq R, \t^{z}_G \le 2 R\bigr)
 +\sum_{u=0}^{R-1} \P\bigl(\t^{z}_G >u\,,\, \t^{z}_{x_0}=u\bigr)
\P\bigl(\t^{x_0}_G \le 2 R-u\bigr) \;.
\end{eqnarray*}
The last inequality results from Markovianness and monotonicity.
Hence, if $z\in B(x_0,r_0)$
\be{eq:ri0}
\P\bigl(\t^{z}_G \le 2 R\bigr)\;\le\; r_0+r+ \P\bigl(\t^{x_0}_G \le 2 R\bigr) \;.
\ee
where we used \eqref{def.bacasym}. As a consequence, using the definition of $c$, if $r_0<c-\bar c$,
 \be{equation5of3.4.bis}
 y_{0} = 1- \P\bigl(\t^{z}_G \leq R\bigr)  \;\geq\; 1 - \P\Bigl(\t^{z}_G \leq 2 R\Bigr)
 \ge 1- c-r_0 \;> 1-\bar c\;.\nonumber
\ee
and the claim holds.  Assume now that the claim is true for $i$, then, by  (\ref{equation2of3.4}) and the inductive hypothesis,
\be{equation4of3.4}
y_{i+1}\;\geq\; 1 -\frac{c}{y_{i}} \;>\;1 -\frac{c}{1-\bar c}\;=\; 1-\bar c.
\ee
The last identity follows from the equality in (\ref{cbarc}).  The claim is proven.
\medskip

To conclude we consider the case in which $k R\leq S \leq (k+1) R$, with $k=\lfloor\frac{S}{R}\rfloor$.  By monotonicity,
$$
\P\bigl(\t^{z}_G > S + R\bigr)\;\geq\; \P\bigl({\t^{z}_G > (k+2) R}\bigr),
\quad \quad\P\bigl(\t^{z}_G >  k R\bigr)\geq \P\bigl({\t^{z}_G > S}\bigr)
$$
Then, by the previous claim and the definition of $\bar c$,
\be{equation5of3.4}
\frac{\P\bigl(\t^{z}_G > S + R\bigr)}{\P\bigl({\t^{z}_G > S}\bigr)}\;\geq\;
\frac{\P\bigl(\t^{z}_G > (k+2) R\bigr)}{\P\bigl({\t^{z}_G > k R}\bigr)}
\;=\;y_{k+1}\,y_k\;>\;(1-\bar c)^2\;=\; 1-\bar c-c\;.
\ee
This conclude the proof of (\ref{eq:stimacebarc}).
\epr

\br{kl2.}
When $c$ is small Lemma \ref{kl4} implies that  the distribution of $\t^{x_0}_G$
does not change very much passing from $S$ to $S+R$.
We have to note that the control here  is given
by estimating near to one the ratios
$\frac{\P\Bigl(\t^{x_0}_G > S + R\Bigr)}{\P\Bigl({\t^{x_0}_G > S}\Bigr)}$ and
$\frac{\P\Bigl(\t^{x_0}_G \le S + R\Bigr)}{\P\Bigl({\t^{x_0}_G \le S}\Bigr)}$,
providing in this way  a ``multiplicative error" on the distribution function. In general such an estimate  is different and more difficult w.r.t.
an estimate with an ``additive error", i.e., which amounts to show that the difference
$|{\P\Bigl(\t^{x_0}_G > S + R\Bigr)}-{\P\Bigl({\t^{x_0}_G > S}\Bigr)}|$ is near
to zero.
The estimate with a multiplicative error is crucial when considering  the tail of the distribution, since  in Lemma \ref{kl4} there are no upper restriction on $S$.
Similarly in (\ref{stimaleq}) the estimate on $\frac{\P\Bigl(\t^{x_0}_G \le S + R\Bigr)}{\P\Bigl({\t^{x_0}_G \le S}\Bigr)}$
is relevant for  $S$ not too small when we can prove, by Lemma \ref{kl3},  that
${c\over \P({\t^{x_0}_G \le S })}$ is small. Note also that the last term in (\ref{stimaleq}) is small due to  Lemma \ref{kl3}.
\er

We use this lemma to prove a multiplicative-error strengthening of Lemma \ref{kl2}.

\bl{aggiuntivo}
Let $(x_0,G)$ be a reference pair satisfying $Rec(R,r)$, $S>R$
and $r$ and $c$ so that
\be{smallness}
r+c+\bar c\;<\;1\;.
\ee
Then, for any $z\in B(x_0,r)$,
\be{improved}
{\P\Bigl({\t^{z}_G > (t+s) S} \Bigr)
\over  \P\Bigl({\t^{z}_G > t S }\Bigr)
  \P\Bigl({\t^{x_0}_G > s  S}\Bigr)}\;\ge\; 1-(c+\bar c+r)\;.
  \ee
  and, for any $z\in\cX$
\be{improved.1}
{\P\Bigl({\t^{z}_G > (t+s) S} \Bigr)
\over  \P\Bigl({\t^{z}_G > t S }\Bigr)
  \P\Bigl({\t^{x_0}_G > s  S}\Bigr)}\;\le\;1+
  {c+\bar c+r\over 1-(c+\bar c+r)}
\ee
\el

\bpr
Inequality \eqref{improved} is a consequence of the top inequality in Lemma \ref{kl2} and inequality \eqref{eq:stimacebarc} of Lemma \ref{kl4}.

To prove \eqref{improved.1} we resort to the decomposition \eqref{t*1tre} which we further decompose by writing
\be{eq:re-further}
\Bigl\{\t^x_G > s  S\,;\,\t^x_{x_0,G}>R\Bigr\}\;=\;\bigcup_{n= 1}^{\floor{sS/R}}
\Bigl\{\t^x_G > s  S\,;\,\t^x_{x_0,G}\in V_n \Bigr\}\;,
\ee
with
\[
V_n\;=\;\left\{\begin{array}{ll}
\bigl( n R,(n+1)R\bigr] & 1\le n \le \floor{sS/R}-1\\
(nR,\infty) & n=\floor{sS/R}
\end{array}\right.
\]
We obtain,
\begin{eqnarray}\label{t*1tre.1}
  \P\Bigl({\t^{z}_G > (t+s)  S}\Bigr)
  &=&\sum_{u=0}^{R}
  \P\Bigl({
      \t^*(t S) = t  S + u \; ;\; \t^{z}_G > t  S + u}\Bigr)\,
  \P\Bigl({\t^{x_0}_G > s  S - u}\Bigr)   \\
   \quad &+&\;
  \sum_{x \in  \{x_0,G\}^c }
  \P\Bigl({\t^{z}_G > t  S        \; ;\; X_{t  S}=x}\Bigr)\,
  \sum_{n= 1}^{\floor{sS/R}}  \P\Bigl({\t^x_G > s  S \; ;\; \t^x_{x_0,G}\in V_n}\Bigr)\nonumber
\end{eqnarray}
and, by monotonicity,
\begin{eqnarray}
\label{eq:rr15}
\P\Bigl(\t^{z}_G > (t+s)  S\Bigr)
&\le&  \P\Bigl(\t^{z}_G > t  S\Bigr)\, \P\Bigl(\t^{x_0}_G > s  S - R \Bigr)\\
  \quad {}&+&\P\Bigl(\t^{z}_G > t S \Bigr)
  \sum_{n= 1}^{\floor{sS/R}}  \sup_{x\in\{x_0,G\}^c}
 \P\Bigl(\t^{x}_G > s  S \,; \; \t^x_{x_0,G}\in V_n\Bigr)\;.\nonumber
\end{eqnarray}
Markovianness and monotonicity imply the following bounds,
\begin{eqnarray}
 \P\Bigl(\t^{x}_G > s  S , \; \t^x_{x_0,G}\in V_n\Bigr)&\le&
\sum_{u\in V_n}  \P\Bigl(\t^{x}_{x_0}=u \, ; \; \t^{x_0}_G \ge sS-u\Bigr)\nonumber\\[5pt]
\le \sum_{u\in V_n}  \P\bigl(\t^{x}_{x_0}=u\bigr)\, \P\bigl(\t^{x_0}_G \ge sS-u\Bigr)
&\le&  \P\Bigl(\t^{x}_{\{x_0,G\}} > nR\Bigr)\, \P\Bigl(\t^{x_0}_G \ge s  S - nR \Bigr)\;.
\end{eqnarray}
Hence, using bounds \eqref{recR+} and (\ref{eq:stimacebarc}),
\be{eq:rr20}
 \P\Bigl(\t^{x}_G > s  S , \; \t^x_{x_0,G}\in V_n\Bigr) \;\le\; r^n\,
 \frac{\P\bigl(\t^{x_0}_G \ge s  S \bigr)}{(1-c-\bar c)^n}\;.
\ee
Replacing this into \eqref{eq:rr15} yields
\begin{eqnarray}
\lefteqn{\P\Bigl(\t^{z}_G > (t+s)  S\Bigr)}\nonumber\\
&\le& \P\bigl({ \t^{z}_G > s  S}\bigr)\,
  \P\bigl({\t^{x_0}_G > t  S - R}\bigr)+\P\bigl(\t^{z}_G > s  S \bigr)\,
\P\bigl(\t^{x_0}_G > t S \bigr)   \sum_{n=1}^\infty \Bigl( {r\over 1-c-\bar c}\Bigr) ^n\nonumber\\
&\le&\P\bigl(\t^{z}_G > s  S \bigr)\, \P\bigl(\t^{x_0}_G > t S \bigr)
\Big[ {1\over 1-c-\bar c}+{a\over 1-a}\Big]
\end{eqnarray}
with $a=r/(1-c-\bar c)$.  [We used (\ref{eq:stimacebarc}) in the first summand.]
Equation (\ref{improved}) follows by noting that
$$
\Big[ {1\over 1-c-\bar c}+{a\over 1-a}\Big]\;=\;1+{c+\bar c\over 1-c-\bar c}+{r\over 1-c-\bar c-r}
\;\le\; 1+{c+\bar c+r\over 1-c-\bar c-r}\;.
$$
\epr

The preceding lemma implies the following improvement of Corollary \ref{cor-fatt}:
\bc{klexp}
Let $(x_0,G)$ be a reference pair satisfying $Rec(R,r)$, and
let $r$ and $c$ be sufficiently small so that $r+c+\bar c<{1\over 2}$ and define
\be{eq:rdelta0}
\delta_0\;:=\;\ln \Big[1+{c+\bar c+r\over 1-(c+\bar c+r)}\Big]\\
\ee
Then, for $S>R$ and $k=1,2,\ldots$
\be{eq:rdiff}
e^{-\d_0 \, k}\;\le\;{\P\Bigl(\t^{x_0}_G > kS \Bigr)\over \P\Bigl(\t^{x_0}_G > S \Bigr)^k}
\;\le\;e^{\d_0 \,k}\;.
\ee
\ec
The inequalities are just an iteration of \eqref{improved} and \eqref{improved.1}. In fact, the left inequality can be extended to any $z\in B(x_0,r)$ and the right inequality to any $z\in\cX$.  We will not need, however, such generality.

The multiplicative bounds of Lemma \ref{aggiuntivo} can be transformed into bounds for $z$ in $ B(x_0,r)$ with the help of the following lemma.

\bl{Claim9}
Consider a reference pair $(x_0, G)$, with $x_0\in\cX$, 
such that $Rec(R,r)$ holds with
$R<T:=\E \t_{G}^{x_0}$.
For all $z\in B(x_0,r_0)$
\be{eq:rclaim8}
 \P\Bigl({\t^{z}_G > tT}\Bigr)\;\ge\;
   \P\Bigl(\t^{x_0}_G > tT\Bigr) \,(1-r-r_0)\;.
   \ee
\el

\bpr
This is proven with a decomposition similar to those used in the proof of Lemma \ref{aggiuntivo}. We have:
\begin{eqnarray*}
 \P\Bigl({\t^{z}_G > tT}\Bigr) &\ge &
 \P\Bigl({\t^{z}_G > tT, \;\t^z_{\{x_0,G\}}<R}\Bigr) = \sum_{u=0}^R \P\Bigl({\t^{z}_G > tT, \;\t^z_{\{x_0\}}=u}\Bigr) \\
 &=& \sum_{u=0}^R \P\Bigl({\t^{z}_G > u, \;\t^z_{\{x_0\}}=u}\Bigr) \,
 \P\Bigl({\t^{x_0}_G > tT-u}\Bigr)\;
\end{eqnarray*}
The last identity is due to Markovianness at $u$. By monotonicity, we conclude
\begin{eqnarray*}
 \P\Bigl({\t^{z}_G > tT}\Bigr) &\ge &
  \P\Bigl({\t^{x_0}_G > tT}\Bigr)\,
  \P\Bigl(\t^{z}_{\{x_0,G\}}\le R, \;\t^z_{\{x_0,G\}}= \t^z_{\{x_0\}}\Bigr) \\
 &=& \P\Bigl(\t^{x_0}_G > tT\Bigr)\,\Bigl[1- \P\Bigl(\bigl\{\t^{z}_{\{x_0,G\}} > R\bigr\}
 \cup\,\bigl\{\t^{z}_{\{x_0,G\}}\not =\t^z_{\{x_0\}}\bigr\}\Bigr) \Bigr]
\end{eqnarray*}
which implies \eqref{eq:rclaim8}.
\epr

The bound \eqref{eq:rclaim8} is complemented by the bound
\be{eq:rclaim8.1}
 \P\bigl({\t^{z}_G > tT}\bigr)\;\le\;
   \P\bigl(\t^{x_0}_G > tT\bigr)\,(1+r)
\ee
implied by the bottom inequality in \eqref{eq:r2}.  Combined with Lemma \ref{aggiuntivo}, the bounds \eqref{eq:rclaim8} and \eqref{eq:rclaim8.1} yield the bound
\be{eq:ex.r1}
\frac{1-(c+\bar c+r)}{1+r}\;\le\;
\frac{\P\Bigl({\t^{z}_G > (t+s) S} \Bigr)}{\P\Bigl({\t^{z}_G > t S }\Bigr)
  \P\Bigl({\t^{z}_G > s  S}\Bigr)}\;\le\;
\frac{1+ \frac{c+\bar c+r}{1-(c+\bar c+r)}}{1-r-r_0}
\ee
valid for all $z\in B(x_0,r_0)$ when $c$, $r$ and $r_0$ are so small that
$c+\bar c+r<1$, $r_0\le r$ and $r+r_0<1$.


\section{Proofs on exponential behaviour}

\subsection{Proof of Theorem \protect\ref{t0}}
We chose an intermediate time scale $S$, with $R< S< T$.  While this scale will finally be related to $\epsilon$ and $r$, for the sake of precision we introduce an additional parameter
\be{eq:rparam}
\eta\;:=\; \frac{S}{T} \quad,\quad 0<\epsilon<\eta<1\;.
\ee
To avoid trivialities in the sequel, we assume that $\eta$, $\epsilon$ and $r$ are small enough so that
\be{eq:rtriv1}
c+\bar c+r \;<\; 1/2 \quad\mbox{and}\quad
\eta+\epsilon+r \;<\;1\;.
\ee
We decompose the proof into nine claims that may be of independent interest.
The claims provide bounds for ratios of the form $\P\bigl(\t^{x_0}_G > \eta \,t \bigr)/e^{-t}$ from which bounds for $\bigl|\P\bigl(\t^{x_0}_G > \eta \,t \bigr)-e^{-t}\bigr|$ can be readily deduced.

\empbf{Claim 1:} Bounds for $t=\eta$.
\be{eq:rclaim1}
e^{-\eta\,\alpha_1}\;\le\;\frac{\P\Bigl(\t^{x_0}_G > \eta \,T \Bigr)}{e^{-\eta}}\;\le\;e^{\eta\,\alpha_0}
\ee
with
\be{eq:r101}
\alpha_0\;:=\; 1 +\frac{1}{\eta}\,\log\Bigl[\frac{1}{1+\eta-2\epsilon} + r\Bigr]\quad\quad \quad
\alpha_1\;:=\; -1 -\frac{1}{\eta} \log\bigl(1-\eta-\epsilon-r\bigr)\;.
\ee
This is just a rewriting of the bounds \eqref{Pt<R} and \eqref{Pt>R} as shown by the following chain of inequalities:
\be{eq:r100}
e^{-\eta[1+\alpha_1]}\;:=\; 1-\eta-\epsilon-r\;\le\;\P\Bigl(\t^{x_0}_G > \eta \,T \Bigr)
\;\le\; 1-\frac{\eta-2\epsilon}{\eta-2\epsilon +1} +r \;=:\; e^{-\eta[1-\alpha_0]}\;.
\ee
Note that, as $\eta$, $\epsilon$ and $r$ tend to zero,
\be{eq:rbcl1}
\alpha_0\,,\,\alpha_1\:=\: O(\epsilon/\eta) + O(r/\eta)\;.
\ee

 \medskip\par\noindent
\empbf{Claim 2:}  Bounds for all $\epsilon<t<\eta$,
\be{eq:rclaim2}
C_-^{(1)}\,e^{-t\,\alpha_1}\;\le\;\frac{\P\Bigl(\t^{x_0}_G > t \,T \Bigr)}{e^{-t}}
\;\le\; C_+^{(1)}\,e^{t\,\alpha_0}
\ee
with functions 
$C_\pm^{(1)}(\eta,\epsilon,r) \;=\; 1+ O(\eta)+O(\epsilon) + O(r)\;$.
Indeed, inequalities \eqref{eq:r100} hold also with $\eta$ replaced by $t$ and yield
\be{eq:rclaim2.1}
e^{-t\alpha_1}\,e^{t[\alpha_1-\alpha_1^t]}\;=\; e^{-t\,\alpha_1^t}\;\le\;\frac{\P\Bigl(\t^{x_0}_G > t \,T \Bigr)}{e^{-t}}
\;\le\; e^{t\,\alpha_0^t}\;=\; e^{t\,\alpha_0}\,e^{t[\alpha_0^t-\alpha_0]}
\ee
where $\alpha_0^t$ and $\alpha_1^t$ are defined as in \eqref{eq:r101} but with $t$ replacing $\eta$.  This is precisely \eqref{eq:rclaim2} with
\be{eq:rcplus1}
C_+^{(1)}\;=\; \sup_{\epsilon<t\le \eta} e^{t[\alpha_0^t-\alpha_0]}
\;=\; \sup_{\epsilon<t\le \eta} \frac{(1+t-2\epsilon)^{-1}+r}{\bigl[(1+\eta-2\epsilon)^{-1}+r\bigr]^{t/\eta}}
\ee
and
\be{eq:rcminus1}
C_-^{(1)}\;=\; \inf_{\epsilon<t\le \eta} e^{t[\alpha_1-\alpha_1^t]}
\;=\; \inf_{\epsilon<t\le \eta} \frac{\bigl[1-\eta-\epsilon-r\bigr]^{t/\eta}}{1-t-\epsilon-r}\;.
\ee

 \medskip\par\noindent
\empbf{Claim 3:}  Bounds for all $t\le \epsilon$,
\be{eq:rclaim3}
C_-^{(0)}\,e^{-t\,\alpha_1}\;\le\;\frac{\P\Bigl(\t^{x_0}_G > t \,T \Bigr)}{e^{-t}}
\;\le\; C_+^{(0)}\,e^{t\,\alpha_0}\;,
\ee
with functions 
$C_\pm^{(0)}(\eta,\epsilon,r) \;=\; 1+ O(\eta)+O(\epsilon) + O(r)\;$.
Indeed, the restriction $t\le \epsilon$ implies that
\be{eq:rclaim3.0}
1-(2\epsilon+r)\;\le\; \P\Bigl(\t^{x_0}_G > \epsilon \,T \Bigr) \;\le\;
\P\Bigl(\t^{x_0}_G > t \,T \Bigr)\;\le\;1\;,
\ee
where the leftmost inequality is a consequence of \eqref{Pt<R} plus the right continuity of probabilities.  The upper bound in \eqref{eq:rclaim3} is a consequence of the rightmost inequality in \eqref{eq:rclaim3.0} and the inequalities
\be{eq:rclaim3.1}
\frac{\P\Bigl(\t^{x_0}_G > t \,T \Bigr)}{e^{-t}}\;\le\;e^t\;\le\; e^\epsilon\;\le\; e^\epsilon
\,e^{t\,\alpha_0}\; =\; \; C_+^{(0)}\,e^{t\,\alpha_0}
\ee
with
$C_+^{(0)}\;=\; e^\epsilon$.
The lower bound in \eqref{eq:rclaim3}, in turns, follows form the leftmost inequality in \eqref{eq:rclaim3.0} and the inequalities
\be{eq:rclaim3.3}
\frac{\P\Bigl(\t^{x_0}_G > t \,T \Bigr)}{e^{-t}}\;\ge\; e^t \,\bigl[1-(2\epsilon+r)\bigr]\;\ge\;
e^{-t\,\alpha_1} \,\bigl[1-(2\epsilon+r)\bigr]\;=\; C_-^{(0)}\, e^{-t\,\alpha_1}
\ee
with
$C_-^{(0)}\;=\; 1-(2\epsilon+r)$.
\medskip\par
Putting the preceding three claims together we readily obtain

 \medskip\par\noindent
\empbf{Claim 4:}  Bounds for all $t\le \eta$,
\be{eq:rclaim4}
\overline C_-\,e^{-t\,\alpha_1}\;\le\;\frac{\P\Bigl(\t^{x_0}_G > t \,T \Bigr)}{e^{-t}}
\;\le\; \overline C_+\,e^{t\,\alpha_0}\;,
\ee
with functions
$\overline C_\pm(\eta,\epsilon,r) \;=\; 1+ O(\eta)+O(\epsilon) + O(r)$.
Indeed, this follows from the three preceding claims, putting
 \be{eq:rclaim4.1}
 \overline C_+ = \max\bigl\{ C_+^{(0)},C_+^{(1)}\bigr\}, \quad\quad \quad
 \overline C_- = \min\bigl\{ C_-^{(0)},C_-^{(1)}\bigr\}\;.
 \ee

 \medskip\par\noindent
\empbf{Claim 5:}  Bounds for $t=k\eta $. Let $\delta_0$ be as in Corollary \ref{klexp}, then for any integer $k\ge 1$,
\be{eq:rclaim5}
e^{-k\,\eta\, \lambda_1}\;\le\;\frac{\P\Bigl(\t^{x_0}_G > k\,\eta\,T \Bigr)}{e^{-k\,\eta}}
\;\le\; e^{k\,\eta\,\,\lambda_0}
\ee
with
\be{eq:rclaim6.01}
 \lambda_0=\alpha_0+\frac{\delta_0}{\eta}\;, \quad\quad \quad
  \lambda_1=\alpha_1+\frac{\delta_0}{\eta}\;.
 \ee

 This result amounts to putting together \eqref{eq:rdiff} and \eqref{eq:rclaim1}.
 Notice that ---by \eqref{eq:defbarc}--\eqref{eq:rr7}--- $c,\,\bar c=O(\epsilon) + O(r)$, hence from \eqref{eq:rbcl1} and the definition \eqref{eq:rdelta0} of $\delta_0$,
\be{eq:rj1}
\lambda_0\,,\, \lambda_1=O\bigl(\epsilon/\eta\bigr)+O\bigl(r/\eta\bigr)\;.
\ee

 \medskip\par\noindent
\empbf{Claim 6:} Bounds for any $t>0$,
\be{eq:rclaim6}
C_-\,e^{-t\,\lambda_1}\;\le\;\frac{\P\Bigl(\t^{x_0}_G > t \,T \Bigr)}{e^{-t}}
\;\le\; C_+\,e^{t\,\lambda_0}
\ee
with
 \be{eq:rclaim6.1}
 C_+=\overline C_+\,\Bigl[1+\frac{c+\bar c +r}{1-(c+\bar c+r)}\Bigr], \quad\quad \quad
 C_-=\overline C_-\,\bigl[1-(c+\bar c+r)\bigr]\;.
 \ee

 Indeed, for any $t>0$ there exist an integer $k\ge 0$ and $0\le t'<\eta$ such that
 $t=k\,\eta + t'$.  Hence,
 \be{eq:rclaim6.3}
 \frac{\P\Bigl(\t^{x_0}_G > t \,T \Bigr)}{e^{-t}}\;=\;
 \frac{\P\Bigl(\t^{x_0}_G > (k\,\eta+t') \,T \Bigr)}{P\Bigl(\t^{x_0}_G > k\,\eta \,T \Bigr)\,
 P\Bigl(\t^{x_0}_G > t' \,T \Bigr)} \,
 \frac{\P\Bigl(\t^{x_0}_G > k\,\eta\,T \Bigr)}{e^{-k\,\eta}}\,
 \frac{\P\Bigl(\t^{x_0}_G > t '\,T \Bigr)}{e^{-t'}}
 \ee
and the claim follows from \eqref{improved}. \eqref{improved.1}, \eqref{eq:rclaim4} and \eqref{eq:rclaim5}.


 \medskip\par\noindent
\empbf{Claim 7:} The bounds \eqref{eq:rclaim6} implies the bound \eqref{expquant} for any $t>0$.

Indeed, subtracting 1 and multiplying through by $e^{-t}$, \eqref{eq:rclaim6} implies
\be{eq:rquant}
\Big|\P\big({\t_{G}^{x_0}} >t\,T\big)-e^{-t}\Big|\;\le\; C \, e^{-t}
\ee
with
\be{eq:rquant1}
C\;=\;  \; \max\Bigl\{\bigl(C_+\,e^{t\,\lambda_0}-1\bigr)\,,\, \bigl|1-C_-\,e^{-t\,\lambda_1}\bigr|\Bigr\}\;.
\ee

To conclude the proof we notice that
$C_+\,,\, C_-= 1+O(\eta)+O(\epsilon)+O(r)$ and, hence, by \eqref{eq:rj1},
%
\be{eq:rquant10}
C\;=\; O(\eta)+O(\epsilon/\eta)+O(r/\eta)\;.
\ee
At this point we choose $\eta$ appropriately to satisfy the asymptotic behavior \eqref{eq:rconv}. A democratic choice, that makes the different contributions of comparable size is
\be{eq:rdif}
\eta\;=\; \sqrt{\max\{\epsilon,r\}}\;.
\ee


 \medskip\par\noindent
\empbf{Claim 8:} Let $r_0$ be such that $r+r_0<1$. Then for any $z\in B(x_0,r_0)$ and any $t>0$,
\be{eq:rclaim7}
\widetilde C_-\,e^{-t\,\lambda_1}\;\le\;\frac{\P\Bigl(\t^{z}_G > t \,T \Bigr)}{e^{-t}}
\;\le\; \widetilde C_+\,e^{t\,\lambda_0}
\ee
with
 \be{eq:rclaim7.1}
 \widetilde C_+=\ C_+\biggl[1+\frac{c+\bar c +r}{1-(c+\bar c+r)}\biggr] \quad \quad \quad
 \widetilde C_-= C_-\,\bigl(1-r-r_0\bigr)\;.
 \ee

Indeed, the upper bound follows from the upper bound in \eqref{eq:rclaim6} and \eqref{improved.1} with the substitutions $S\to T$, $t=0$ and $s\to t$.  The lower bound is a consequence of the lower bound in \eqref{eq:rclaim6} and \eqref{eq:rclaim8} in Lemma \ref{Claim9}
This concludes also the proof of \eqref{ct0}.
\eproof



\subsection{Proof of Theorem \protect\ref{t1}:}

\paragraph{Part (I).}  The results follow from \eqref{expquant} and \eqref{ct0} by letting
both $\e_n:={R_n\, /\, \E \t_{G^{(n)}}^{(n),x_0^{(n)}} }$ and $r_n$ tend to zero.

\paragraph{Part (II) (i).}
Let
\be{eq:rn100}
\vartheta_n\;:=\; \frac{\t_{G^{(n)}}^{(n),x_0^{(n)}}}{Q_n(\z)}\;.
\ee
Combining the definition \eqref{defQzeta} of $Q_n(\z)$ with Corollary \ref{klexp} we obtain
\be{eq:rn101}
\bigl[\z\,e^{-\delta_{0,n}}\bigr]^k\;\le\; \P\bigl(\vartheta_n >k\bigr)\;\le\;
\bigl[\z\,e^{\delta_{0,n}}\bigr]^k
\ee
with
\be{eq:rn102}
\delta_{0,n} \longrightarrow 0 \quad \mbox{as} \quad \epsilon_n, r_n \to 0\;.
\ee
A simple argument based on Markov inequality [see \eqref{eq:r200}--\eqref{eq:r2001.1}   below] shows that hypothesis Hp.$G(Q_n(\zeta))$ implies hypothesis Hp.$G(T_n)$ [in fact, they are equivalent, as shown in Theorem \ref{t3}].  Hence, \eqref{eq:rn102} holds under hypothesis Hp.$G(Q_n(\zeta))$ and \eqref{eq:rn101} shows that the sequence of random variables $(\vartheta_n)$ is
exponentially tight.  Therefore, the sequence is relatively compact in the weak topology (Dunford-Pettis theorem) and every subsequence has a sub-subsequence that converges in law.  Let
$(\vartheta_{n_k})_k$ be one of these convergent sequences and let
$\vartheta$ be its limit.
Taking limit in \eqref{eq:r2} we see that
\be{preexp}
\P( \vartheta > t+s) \;=\; \P( \vartheta > t)\,\P( \vartheta > s)
\ee
for all continuity points $s,t > 0$.
  Since these points are dense
  and the distribution function is right-continuous,
   \eqref{preexp} holds
  for all $s,t\ge 0$.  Furthermore, the limit of \eqref{eq:rn101} implies that
  \be{eq:rn103}
\P\bigl(\vartheta>k\bigr)\;=\; \z^k\;.
\ee
 We conclude that $\vartheta$ is an exponential variable of rate $-\log\z$.  As every subsequence of  $(\theta_n)$ converges to the same $\exp(-\log\z)$ law, the whole sequence does.

 \paragraph{ Part (II) (ii).}
 The rightmost inequality in \eqref{eq:rn101} implies that fixing $\z<\widetilde\z<1$, we have that  for $n$ large enough
 \[
 \P\bigl(\vartheta_n >k\bigr)\;\le\; \widetilde\z^k
 \]
 for all $k\ge 1$.  By monotonicity this implies that
 \be{eq:rn104}
 \P\bigl(\vartheta_n >t\bigr)\;\le\; \widetilde\z^{\floor{t}}
 \ee
for $n$ large enough.  As the function in the right-hand side is integrable, we can apply dominated convergence and conclude:
 \be{eq:rn105}
  \lim_{n \to \infty}{Q_n(\zeta)\over T_n^E }\;=\;
  \lim_{n\to\infty} \int_0^\infty  \P\bigl(\vartheta_n >t\bigr)\,dt \;=\;
  \int_0^\infty  \P\bigl(\vartheta>t\bigr)\,dt \;=\; -\ln\z\;.\qquad\qed
 \ee

\subsection{Proof of Theorem \protect\ref{thEEquant}}

We assume $\epsilon$ and $r$ small enough so that $c+\bar c+ r<1/2$.
We shall produce an upper and a lower bound for the ratio
\be{eq:rn106}
REE(S)\;:=\;{\P\Bigl({\t^{x_0}_G}\in (kS,(k+1)S)\Bigr)\over
\P\Bigl({\t^{x_0}_G} > S \Bigr)^k \P\Bigl({\t^{x_0}_G} \le S \Bigr)}
\ee
for $k<T/S$.
We start with the identity
\begin{eqnarray}
\label{eq:rn107}
\P\Bigl({\t^{x_0}_G}\in (kS,(k+1)S)\Bigr) &=&
\P\Bigl({\t^{x_0}_G}\in (kS,(k+1)S), \t^*(kS-R)\le kS\Bigr)\nonumber\\
&&\; {}+ \P\Bigl({\t^{x_0}_G}\in (kS,(k+1)S), \t^*(kS-R)> kS\Bigr)\quad
\end{eqnarray}

Applying Markov property and monotonicity we obtain the upper bound
\begin{eqnarray}
\label{eq:rn108}
\P\Bigl({\t^{x_0}_G}\in (kS,(k+1)S)\Bigr) &\le&
\P\Bigl({\t^{x_0}_G}>kS-R\Bigr)\, \P\Bigl({\t^{x_0}_G}\le R+S\Bigr)\nonumber\\
&& \quad{}+ \P\Bigl({\t^{x_0}_G}>kS-R\Bigr) \,r\quad
\end{eqnarray}
and the lower bound
\be{eq:rn109}
\P\Bigl({\t^{x_0}_G}\in (kS,(k+1)S)\Bigr)
\;\ge\; \P\Bigl({\t^{x_0}_G}>kS\Bigr)\, \P\Bigl({\t^{x_0}_G}\le S\Bigr)-
\P\Bigl({\t^{x_0}_G}>kS\Bigr)\, r\;.
\ee
This lower bound and the leftmost inequality in \eqref{eq:rdiff} yields the lower bound
\begin{eqnarray}
\label{eq:rn110}
REE(S)&\ge&
{\P\Bigl({\t^{x_0}_G}>kS\Bigr)\P\Bigl({\t^{x_0}_G}\le S\Bigr) \biggl[1- {r\over \P\Bigl({\t^{x_0}_G}\le S\Bigr)}\biggr]
\over \P\Bigl({\t^{x_0}_G} > S \Bigr)^k \P\Bigl({\t^{x_0}_G} \le S \Bigr)}\nonumber\\[10pt]
&\ge& e^{-\d_0 k}
\biggl[1- {r\over \P\Bigl({\t^{x_0}_G}\le S\Bigr)}\biggr]\;.
\end{eqnarray}

To obtain a bound in the opposite direction we use the upper bound \eqref{eq:rn108}:
\be{eq:rn111}
REE(S)\;\le\; {\P\Bigl({\t^{x_0}_G}>kS-R\Bigr) \Big[\P\Bigl({\t^{x_0}_G}\le R+S\Bigr)+ r\Big]
\over \P\Bigl({\t^{x_0}_G} > S \Bigr)^k \,\P\Bigl({\t^{x_0}_G} \le S \Bigr)}\;.
\ee
We bound the right-hand side through the rightmost inequality in \eqref{eq:rdiff} and the following two easy consequences of \eqref{eq:stimacebarc}:
\be{eq:rn112}
\P\Bigl({\t^{x_0}_G}>kS-R\Bigr) \;\le\; \frac{\P\Bigl({\t^{x_0}_G}>kS\Bigr)}{1-c-\bar c}
\ee
[obtained by replacing $S \to kS-R$ in \eqref{eq:stimacebarc}] and
\be{eq:rn113}
\P\Bigl({\t^{x_0}_G}\le R+S\Bigr)\;\le\;
\P\Bigl({\t^{x_0}_G}\le S\Bigr)(1+c+\bar c)\;.
\ee
The result is
\be{eq:rn114}
REE(S)\;\le\; {1\over 1-c-\bar c}e^{\d_0 k}\biggl[1+c+\bar c+{r\over \P\Bigl({\t^{x_0}_G}\le S\Bigr)}\biggr]
\ee

Inspection shows that the dominant order in both bounds \eqref{eq:rn110} and \eqref{eq:rn114}
is given by the term
\be{eq:rn115}
{r\over \P\Bigl({\t^{x_0}_G}\le S\Bigr)}\;\le\; \frac{r}{1-e^{-\eta(1-\alpha_0)}}\;,
\ee
the last inequality being a consequence of the rightmost inequality in \eqref{eq:rclaim1}.  From \eqref{eq:r101} we see that
\be{eq:rn116}
\frac{r}{1-e^{-\eta(1-\alpha_0)}} \;=\; \frac{r}{O(\eta)+O(\epsilon)+O(r)} \;=\;O\bigl(\epsilon/\eta\bigr) +
O\bigl(r/\eta\bigr)\;.
\ee
With this observation, the bounds \eqref{eq:rn110} and \eqref{eq:rn114} imply
\be{eq:rn117}
\bigl|REE(S)-1\bigr| \;\le\; O\bigl(\epsilon/\eta\bigr) +
O\bigl(r/\eta\bigr)\;. \qquad\qed
\ee

\section{Proof of the relation between different hypotheses for exponential behavior}

In subsection \ref{proofTcomparison} we prove Theorem \ref{t3} and in subsection  \ref{abc} we give examples that show that the converse of each of the first two implications in \eqref{eq:implic} are false.

\subsection{Proof of Theorem \ref{t3}}\label{proofTcomparison}
\medskip

\noindent
\emph{(i) } $Hp.A\Rightarrow Hp.G^{LT}$:
By Markovianness, if $x\not= y\in\cX$,
$$
\P(\xi^x_y(x)>t)\;=\;\P(\xi^x_y(x)>t-1)\,\P(\xi^x_y(x)\ge 1)\;=\;\P(\xi^x_y(x)>t-1)\,\P(\tilde \t^x_y>\tilde \t^x_x)\;,
$$
thus
$$
\P(\xi^x_y(x)>t)\;=\;\P(\tilde \t^x_y>\tilde \t^x_x)^t
$$
and
$$
\E(\xi^x_y(x))\;=\;\P(\tilde \t^x_y<\tilde \t^x_x)^{-1}.
$$
Therefore,
\be{rhoAnsecond}\r_A(n)\;=\; {\sup_{z\in\{x_0,G\}^c}\E\Bigl(\xi^z_{\{x_0,G\}}(z)\Bigr)\over \E\Bigl(\xi^{x_0}_G (x_0)\Bigr)}
\;=\;  {\sup_{z\in\{x_0,G\}^c}\E\Bigl(\xi^z_{\{x_0,G\}}(z)\Bigr)\over T_n^{LT}}\;.
\ee
Furthermore, by Markov inequality
\begin{eqnarray*}
\P\Bigl(\t^z_{\{x_0,G\}}>R_n\Bigr)&=&\P\Bigl(\sum_{x\in\{x_0,G\}^c}\xi^z_{\{x_0,G\}}(x)>R_n\Bigr)\\
&\le&
{\sum_{x\in\{x_0,G\}^c}\E\Bigl(\xi^x_{\{x_0,G\}}(x)\Bigr)\over R_n}
\le \frac{|\cX|\,\rho_A(n)\,T_n^{LT}}{R_n}\;.
\end{eqnarray*}
The proposed implication follows, for instance, by choosing $R_n=\e_nT_n^{LT}$ with $\e_n={\sqrt{ |\cX|\rho_A(n)}}$.
\medskip

\noindent
\emph{(ii) } $Hp.G^{LT}\Rightarrow Hp.G^{E}$:
It is an immediate consequence of the obvious inequality $T_n^E>T_n^{LT}$.
\medskip

\noindent
\emph{(iii) } $Hp.G^{E}\Leftarrow   Hp.G^{Q(\zeta)}$:
By Markov inequality,
\be{eq:r200}
\P\Bigl(\t^{x_0}_G>t\Bigr)\;\le\; \frac{\E\bigl(\t^{x_0}_G\bigr)}{t}\;.
\ee
Thus,
\be{eq:r2001.1}
T_n^{Q(\zeta)}\;\le\; \zeta^{-1}\,T_n^E\;.
\ee
\medskip

\noindent
\emph{(iv) } $Hp.G^{E}\Rightarrow   Hp.G^{Q(\zeta)}$:
We bound:
\begin{equation}\label{eq:r201}
\E\Bigl(\t^{x_0}_G\Bigr) \;=\; \sum_{k\ge 0} \quad \sum_{kT_n^{Q^\zeta}\le t< (k+1)T_n^{Q^\zeta}} \quad
\P\Bigl(\t^{x_0}_G>t\Bigr) \nonumber\\
\;\le \; T_n^{Q(\zeta)} \sum_{k\ge 0} \P\Bigl(\t^{x_0}_G>k\,T_n^{Q(\zeta)} \Bigr)\;.
\end{equation}
Hence, by Lemma \ref{klexp},
\be{eq:r202}
T_n^E\;\le\, T_n^{Q(\zeta)} \Bigl[ 1 + O\Bigl((\zeta+\e_n+r_n)\Bigr)\Bigr]\;.
\ee

\medskip

\noindent
\emph{(v) } $Hp.G^{E}  \Leftarrow Hp.B$:
Apply the Markov inequality to $\P(\t^x_{\{{x_0},G\}}>R_n)$ and choose $R_n={\sqrt{T_n^E\sup_{z\in\{{x_0},G\}^c}\E(\t^z_{\{{x_0},G\}}})}$.
\medskip

\noindent
\emph{(vi) } $Hp.G^{E}  \Rightarrow Hp.B$:
\begin{eqnarray*}
\sup_{z\in\{{x_0},G\}^c}\E(\t^z_{\{{x_0},G\}})= \sup_{z\in\{{x_0},G\}^c}\sum_{t\ge 0}\P(\t^z_{\{{x_0},G\}}>t)
\le R_n+R_n \sup_{z\in\{{x_0},G\}^c}\sum_{N=1}^\infty\P(\t^z_{\{{x_0},G\}}>NR_n)\\
\le R_n+R_n \sum_{N=1}^\infty\sup_{z\in\{{x_0},G\}^c}\P(\t^z_{\{{x_0},G\}}>R_n)^N
\le R_n+{R_n\over 1-r_n}
\le const. \, R_n
\end{eqnarray*}
As $Hp.G^{E} $ implies $R_n\prec T^E_n$ we obtain that $\rho_B(n)\to 0$.

\subsection{Counterexamples}\label{controesempi}
\subsubsection{h model (freedom of starting point)\label{hmodel}}

We show that even in the finite--volume
case  the set of starting points that verify Hp. B is in general bigger
than the set of good starting points for Hp. A and that these are
not limited to the deepest local minima of the energy.

\begin{figure}
\centering
\includegraphics[width=0.5\linewidth]{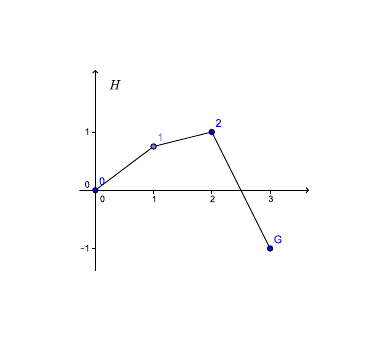}
\caption{The h model.}
\label{fig:hmodel}
\end{figure}

Let us consider a birth-and-death model with state space
$\{0,1,2,G\}$ and  transition matrix
\begin{align}
P(0,1)&:=\frac{1}{2}p^{h}&  P(1,2)&:=\frac{1}{2}p^{1-h}&  P(1,0)&:=\frac{1}{2}&\nonumber \\
P(2,G)&:=\frac{1}{2}&       P(2,1)&:=\frac{1}{2}&         P(G,2)&:=\frac{1}{2}p^{2}&\label{Ph}
\end{align}
and $P(i,i)=1-\sum_{j\not=i}P(i,j)$ for $i\in\{0,1,2,G\}.$
This model can be seen as a Metropolis chain with energy function $H(0)=0,\ H(1)=h,\ H(2)=1,\ H(G)=-1$
(see fig. \ref{fig:hmodel}).

By Lemma 3.3 in \cite{BM}, the reference pair $(0,G)$ always fulfils
Hp. A (and thus Hp. B), since $0$ is the most stable local minimum of the
energy.
We focus here on the pair $(1,G)$ and show that for some values of the parameter $h$ this pair verifies Hp. B or even Hp. A although $1$ is not a local minimum of the energy.

It is easy to show (see (\ref{xih}), (\ref{th1}), (\ref{th2}) for exact computations) that
\begin{align}
\E(\tau_{\left\{ 1,G\right\} }^{2})&\asymp 1&  \E(\xi_{\left\{ 1,G\right\} }^{2})&\asymp 1&    \E(\tau_{G}^{1})\asymp p^{-1}\nonumber \\
\E(\xi_{G}^{1})&\asymp p^{-1+h}&           \E(\tau_{\left\{ 1,G\right\} }^{0})&\asymp p^{-h} &    \E(\xi_{\left\{ 1,G\right\} }^{0})\asymp p^{-h},
\label{P-1}
\end{align}
where the notation $f\asymp g$ means that for for $p$ small the
ratio $f/g$ is bounded from above and from below by two positive constants.

By (\ref{P-1}),
\[
\rho_{A}=\frac{E(\xi_{\left\{ 1,G\right\} }^{0})}{E(\xi_{G}^{1})}\asymp p^{1-2h};\ \ \rho_{B}=\frac{E(\tau_{\left\{ 1,G\right\} }^{0})}{E(\tau_{G}^{1})}\asymp p^{1-h},
\]
so that Hp. A holds when $h<\frac{1}{2}$ while Hp. B holds for $h<1$.
It is possible to show that in this example Hp. B is optimal, since $h<1$ is also a necessary condition to have the exponential law for
$\t^1_G / \E(\t^1_G)$ in the limit $\b \to \infty$.


\subsubsection{abc model (strength of recurrence properties)\label{abc}}

When the cardinality of the configuration space diverges, condition Hp. $A$ is in general stronger than Hp. $G^{LT}$ which is stronger than Hp. $B$.

We show this fact with the help of a simple one-dimensional model, that provides the counterexamples of the
missing implications in theorem \ref{t3} (see also fig. \ref{fig:Venn2}):

We consider
a class of birth and death models that we use to discuss the differences
among the hypotheses. The models are characterized by three positive
parameters $a$, $b$, $c$.

Let $\Gamma:=\left\{ 0,\ldots,L\right\} $ be the state space. We take $L=n\to \infty$.
Let
$a,b,c$ be real parameters, with $b\le a$, $b\le c$.

The nearest-neighbor transition probabilities are defined as:

\begin{flalign}
P_{0,1}: & =\frac{1}{2}L^{-a}\nonumber \\
P_{x,x+1}: & =\frac{1}{2}\text{ for }x\in[1,\, L-3] &  & P_{x,x-1}:=\frac{1}{2}\text{ for }x\in[1,\, L-3]\nonumber \\
P_{L-2,L-1}: & =\frac{1}{2}L^{-c} &  & P_{L-2,L-3}:=\frac{1}{2}L^{-b}\label{P}\\
P_{L-1,L}: & =\frac{1}{2} &  & P_{L-1,L-2}:=\frac{1}{2}\nonumber \\
 &  &  & P_{L,L-1}:=\frac{1}{2}L^{-2(a+b+c)}\nonumber
\end{flalign}
 with $P_{x,x}:=1-P_{x,x-1}-P_{x,x+1}$and $P_{x,y}=0$ if $|x-y|>1$.

\begin{figure}
\centering
\includegraphics[width=0.7\linewidth]{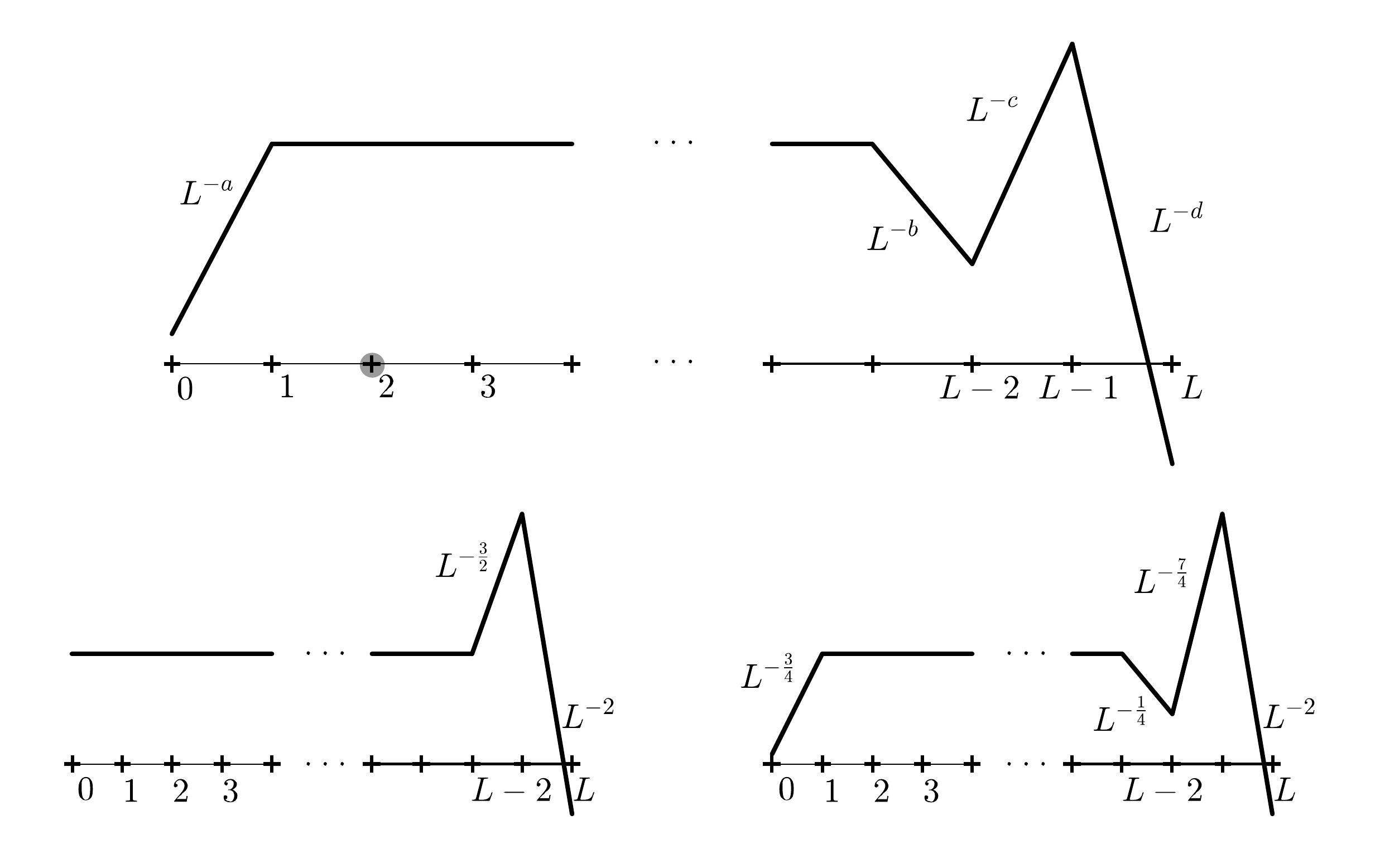}
\caption{Transition dyagrams for the {\it abc} model in the general case and for the choices $(0,0,\frac 3 2)$ and $(\frac 3 4, \frac 1 4, \frac 7 4)$.}
\label{fig:abc2}
\end{figure}
The unique equilibrium measure  can be obtained
trough reversibility condition:

\begin{align}
\mu(0) &=L^{a}/\mathit{\mathcal{Z}}&
\mu(x) &=1/\mathit{\mathcal{Z}}\; \forall  x\in[1,\, L-3]&  &\quad\quad&    \label{eq:pi}\\
\mu(L-2) &=L^{b}/\mathit{\mathcal{Z}}&
\mu(L-1) &=L^{b-c}/\mathit{\mathcal{Z}}&
\mu(L) &=L^{2a+3b+c}/\mathit{\mathcal{Z}},& \nonumber
\end{align}
where $\mathit{\mathcal{Z}}$ is the normalization
factor.

\bigskip\noindent
{\Large\bf Example 1 ($\mathbf{\frac 5 8,\frac{1}{4},\frac{7}{4}}$ model )}

%
In this example, for the reference pair $(0,L)$, $Hp. G^{LT}$ holds whereas $Hp. A$ does not.

The state space is $[0,L]$, with all birth and death probabilities
equal to $\frac{1}{2}$ except for $P_{0,1}:=\frac{1}{2}L^{-\frac{5}{8}}$,
$P_{L-2,L-3}:=\frac{1}{2}L^{-\frac{1}{4}}$, $P_{L-2,L-1}:=\frac{1}{2}L^{-\frac{7}{4}}$
and $P_{L,L-1}:=\frac{1}{2}L^{-\frac{21}{2}}$.

For rigorous computations see the appendix, here we discuss the model
at heuristical level.

In view of Markov inequality, in order to check
$Hp. G^{LT}$ it is sufficient to show that
\be{GLTsuff}
\sup_{x\in\{1,L-1\}} \frac{\E (\tau^x_{0,L})}
{\E (\xi^0_{L})}
\longrightarrow 0.
\ee

Since the well in $L-2$ is rather small, the maximum of the mean
resistance times can be estimated as the time needed to cross the
plateau (of order $L^{2}$):
\be{m1xt}
\max_{x\in\{1,L-1\}}\mathbb{E}(\tau_{0,L}^{x})=\mathbb{E}(\tau_{0,L}^{L-2})\asymp L^{2}
\ee
(see \eqref{DAFARE} for a rigorous computation).

The mean local time can be estimated by using \eqref{MLT}:
the probability that the process starting from $0$ visits $L$ before
returning to $0$ can be estimated as the probability of stepping
out from $0$ (which is $P_{0,1}=\frac{1}{2}L^{-\frac{5}{8}}$) times
the probability of reaching $L-2$ before returning (which is of order
$L^{-1})$ times the probability that the process starting from $L-2$
visits $L$ before $0$ (which is of order $L^{\frac{5}{4}}/L^{\frac{7}{4}}$).
The mean local time of $0$ (i.e. the time spent in $0$ before the
visit to $L$) is therefore
\begin{eqnarray}\label{m1lt}
\mathbb{E}(\xi_{L}^{0})
&=&
\mathbb{P}\left(\tilde \tau_{L}^{0}< \tilde \tau_{0}^{0}\right)^{-1} \nonumber\\
&=&
\left(P_{0,1}\mathbb{P}\left(\tilde \tau_{L-2}^{1}<\tilde \tau_{0}^{1}\right)\mathbb{P}\left(\tilde \tau_{L}^{L-2}<\tilde \tau_{0}^{L-2}\right)\right)^{-1}\asymp
L^{\frac 5 8 + 1 +\frac 1 2} =
 L^{\frac{17}{8}}
\end{eqnarray}
(as proven in \eqref{MLT}).
Thus, by \eqref{m1xt}, \eqref{m1lt} the ratio in \eqref{GLTsuff} goes like
$\asymp L^{-\frac 1 8}$ and condition $Hp. G^{LT}$
holds.

On the other hand,
the maximum of the mean local time before hitting $0$ or $L$ is
reached in $L-2$ and it is of the order of the inverse of $P_{L-2,L-3}$
(i.e. the probability of exiting from $L-2$) times the probabily of exiting
the plateau in $0$ (of order $L^{-1})$:
\be{m1xlt}
\max_{x\in[1,L-1]}\mathbb{E}(\xi_{0,L}^{x})=\mathbb{E}(\xi_{0,L}^{L-2})\asymp\left(P_{L-2,L-3}\mathbb{P}\left(\tilde \tau_{0}^{L-3}<\tilde \tau_{L-2}^{L-3}\right)\right)\asymp L^{\frac{5}{4}}
\ee
(see \eqref{eq:stimaLT} for particulars).
Thus, by \eqref{m1lt}, \eqref{m1xlt},
\be{notA}
L \rho_A(L)
\asymp
L \frac {L^{\frac 5 4}}{L^{\frac {17} 8}}
=L^{\frac 1 8}
\longrightarrow \infty
\ee
that is, condition $A$ does not hold.

{\Large\bf Example 2 ($\mathbf{0,0,\frac{3}{2}}$ model )}

%

This choice of parameters corresponds to a birth-and-death chain where all birth
and death probabilities are set equal to $\frac{1}{2}$ except for
$P_{L-2,L-1}:=\frac{1}{2}L^{-\frac{3}{2}}$ and $P_{L,L-1}:=\frac{1}{2}L^{-3}$.

We show that, for the reference pair $(0,L)$, condition $Hp. B$ holds whereas
condition $Hp. G^{LT}$ does not.

Precise computations can be found in the appendix, here we discuss
the model at heuristic level.

We start by estimating $\r_B(L):=\sup_{z<L}{\E ( \t^{z}_{0, L}) / \E (\t^{ 0}_{L}}) $:

Each time the process is in $L-2$, it
has a probability $\frac{1}{2}L^{-\frac{3}{2}}$ of reaching $L$ in two
steps, but the probability to find the process in $L-2$ before the
transition is of order $L^{-1}$; hence,

\[
\mathbb{E}(\tau_{L}^{0})\asymp L\times L^{\frac{3}{2}},\]
see eq. \eqref{eq:T-1} for the exact computation.

The maximum of the mean   times is dominated by a diffusive
contribution:
\[
\max_{x\in[1,L-1]}\mathbb{E}(\tau_{0,L}^{x})=\mathbb{E}(\tau_{0,L}^{L-2})\asymp L^{2},\]
as confirmed by \eqref{DAFARE}.

Thus, $\rho_{B}\asymp L^{-\frac{1}{2}}$and condition $B$ holds.

Next, we show that
\be{notLT}
\P \left( \t^{L/2}_{\{0,L\}} > \E (\xi^0_L) \right) \longrightarrow 1,
\ee
and, therefore, that $Hp. G^{LT}$ does not hold.

Equation \eqref{MLT} allows to compute the mean local time:
The probability that the process starting from $0$ visits $L$ before
returning to $0$ can be estimated as the probability of reaching
$L-2$ before returning (which is of order $L^{-1})$ times the probability
that the process starting from $L-2$ visits $L$ before $0$ (which
is of order $L/L^{\frac{3}{2}}$). The mean local time of $0$ (i.e.
the time spent in $0$ before the visit to $L$) is therefore \[
\mathbb{E}(\xi_{L}^{0})=\mathbb{P}\left(\tilde\tau_{L}^{0}<\tilde\tau_{0}^{0}\right)^{-1}=\left(\mathbb{P}\left(\tilde\tau_{L-2}^{0}<\tilde\tau_{0}^{0}\right)\mathbb{P}\left(\tilde\tau_{L}^{L-2}<\tilde\tau_{0}^{L-2}\right)\right)^{-1}\asymp L^{\frac{3}{2}}\]
(see \eqref{MLT} for a rigorous derivation).

Let $M(t):= \max_{s\le t} \left| \frac L 2 - X_s^{L/2}
\right|$.
Since $x=L/2$ is in the middle of the plateau, we can use
the diffusive bound $M(t) \asymp \sqrt{t}$ for small $t$. By Markov inequality:
$$
\P \left( \t^{L/2}_{\{0,L\}} < L^{\frac 3 2} \right) <
\P \left( M(L^{\frac 3 2}) > \frac L 3 \right)
\le 3 \frac {\E (M(L^{\frac 3 2}))} L \asymp L^{-\frac 1 4 },
$$
that implies \eqref{notLT}.

\appendix
\section{Appendix}
\subsection{Electric networks}

A convenient language to describe the behavior of local and hitting
times in the reversible case exploits the analogy with electric networks.
Here we recall some useful relation between reversible Markov processes
and electric networks. For a more complete discussion, see {[}DS{]}
and references therein.

We associate with a given reversible Markov chain with transition
matrix $P$ a resistance network in the following way:

We call
We call ``resistance'' of an edge $(x,y)$
of the graph associated with the Markov kernel the quantity
\begin{equation}
r_{x,y}:=\left(\mu(x)P_{x,y}\right)^{-1}.\label{defr}
\end{equation}
 Given two disjoint subsets $A,\, B\;\subset\Gamma$, we denote by
the capital letter $R_{B}^{A}$ the total resistance between $A$
and $B$, namely, the total electric current that flows in the network
if we put all the points in $A$ to the voltage $1$ and all the points
in $B$ to the voltage $0$.

It is well-known that the total resistance is related with the mean
local time (see def. \ref{loctime}) of the Markov chain, i.e. with the Green function, by
\begin{equation}
R_{B}^{x}=\frac{\mathbb{E}(\xi_{B}^{x})}{\mu(x)}.\label{RLT}
\end{equation}

Resistances are reversible objects, i.e. $R_{y}^{x}=R_{x}^{y}.$

The voltage at point $y\in\Gamma$ has a probabilistic interpretation
given by
\[
V_{B}^{x}(y)=\mathbb{P}\left(\tilde \tau_{x}^{y}<\tilde \tau_{B}^{y}\right).
\]

\subsection{h model}

\subsubsection{Computation of resistances}

By reversibility, from (\ref{Ph}), we easily get
\be{gibbsh}
\mu(0)=\frac{1}{Z}\qquad
\mu(1)=\frac{p^h}{Z}\qquad
\mu(2)=\frac{p}{Z}\qquad
\mu(G)=\frac{p^{-1}}{Z}\qquad
\ee

Where $Z$ is the renormalization factor. (\ref{gibbsh}) can be seen as the Gibbs measure of the system once taken $\b=-\log p$.

By (\ref{defr}), (\ref{Ph}) and (\ref{gibbsh}), we get
\begin{equation}
r_{k}:=r_{k,k+1}=\frac{\mathcal{Z}}{2}\times\begin{cases}
p^{-h} & \text{ for }k=0\\
p^{-1} & \text{ for }k=1,2
\end{cases}\label{eq:rh}
\end{equation}

\subsubsection{Local times}

In order to compute $\rho_{A}$, we need to estimate the local time
spent in the metastable point $1$ before the transition to $G$ and
the maximum among the local times of the points  $\left\{ 0,2\right\} $
before the transition to $1$ or to $G$.

By (\ref{RLT}), the computation of mean local times is the analogous
of the computation of the total resistance of a series of resistances:
\begin{eqnarray}
\E\left( \xi^0_G \right)&=&\E\left( \xi^2_{1,G} \right)=\mu(0)(r_0+r_1+r_2)\asymp p^{-1}\nonumber\\
\E\left( \xi^1_G \right)&=&\mu(1)(r_1+r_2)\asymp p^{-1+h}\nonumber\\
\E\left( \xi^2_G \right)&=&\mu(0)r_2 \asymp 1 \nonumber\\
\E\left( \xi^2_{1,G} \right)&=&\mu(2)\frac{r_1 r_2}{r_1+r_2}\asymp 1
\label{xih}
\end{eqnarray}

\subsubsection{Hitting times}

Obviously,
\begin{eqnarray}
\E (\t{^0_{1,G}})=\E (\xi^0_1)\asymp p^{-1+h}\nonumber\\
\E (\t{^2_{1,G}})=\E (\xi^2_{,G})\asymp 1\label{th1}.
\end{eqnarray}

More generally, local times provide a useful language to describe the model. E.g.
a relation between hitting times and mean times is
\begin{align}
\mathbb{E}\left(\tau_{A}^{x}\right) & =\mathbb{E}\left(\sum_{t=0}^{\tau_{A}^{x}}\ \sum_{y\notin A}1_{X_{t}^{x}=y}\right)\nonumber \\
 & =\sum_{y\notin A}\mathbb{E}\left(\xi_{A}^{x}(y)\right)=\sum_{y\notin A}\mathbb{E}\left(\xi_{A}^{y}\right)\mathbb{P}\left(\tilde \tau_{y}^{x}< \tilde \tau_{A}^{x}\right),\label{eq:Txi}
\end{align}
where we used the strong Markov property at time $\tilde \tau_{k}^{0}$ in
the last equality. In words, the hitting time is the sum of all local
times of the points visited.

Since $\mathbb{P}\left(\tilde \tau_{1}^{0}<\tilde \tau_{1}^{G}\right)\ge P(1,0)=\frac 1 2 \asymp 1$, by (\ref{eq:Txi}), (\ref{xih}), we see
that

\be{th2}
\mathbb{E}(\tau_{1}^{G})\asymp \mathbb{E}\left(\xi_{0}^{G}\right)\asymp p^{-1}.
\ee
\subsection{abc model}
\subsubsection{Computation of resistances}

By (\ref{defr}), (\ref{P}) and (\ref{eq:pi}), we get for $k\in\left\{ 0,\ldots,L-1\right\} $
\begin{equation}
r_{k}:=r_{k,k+1}=\mathcal{Z}\times\begin{cases}
2 & \text{ for }k\le L-3\\
2L^{c-b} & \text{ for }k>L-3
\end{cases}\label{eq:r}
\end{equation}

Interesting resistances in our one-dimensional model are the total
resistance $R_{0}^{x}$ between $x$ and $0$ and the total resistance
$R_{L}^{x}$ between $x$ and $L$.

\begin{equation}
R_{0}^{x}:=\sum_{k=0}^{x-1}r_{k}=\mathcal{Z}\times\begin{cases}
2n & \text{ for }x\le L-2\\
2(L-2)+2L^{c-b} & \text{ for }x=L-1
\end{cases}.\label{eq:R0n}
\end{equation}

\begin{equation}
R_{L}^{x}:=\sum_{k=x}^{L-1}r_{k}=\mathcal{Z}\times\begin{cases}
2(L-2-x)+4L^{c-b} & \text{ for }x\le L-2\\
2L^{c-b} & \text{ for }x=L-1
\end{cases}.\label{eq:RnL}
\end{equation}

\subsubsection{Local times}

In order to compute $\rho_{A}$, we need to estimate the local time
spent in the metastable point $0$ before the transition to $L$ and
the maximum among the local times of the points in $\left\{ 1,\ldots,L-1\right\} $
before the transition to $0$ or to $L$.

By (\ref{RLT}), the computation of mean local times is the analogous
of the computation of the total resistance of a series of resistances.
By (\ref{eq:RnL}) we get

\be{MLT}
\mathbb{E}(\xi_{L}^{0})=  \mu(0)R_{L}^{0}=2L^{a}\left((L-2)+2L^{c-b}\right)
  \asymp L^{a}\left(L+L^{c-b}\right)(1+o(1)).
\ee

Then, we are interested in the maximum of the local times $\mathbb{E}(\xi_{0,L}^{x})=\mu(x)$$R_{0,L}^{x}$.
Depending on the equilibrium measure, we can have the maximum in the
plateau $[1,\, L-3]$, in the well $L-2$ or (in principle) in the
peak $L-1.$ Since the parallel between two resistances $r$ and $R$,
with $r\le R$ is between $r/2$ and $r$, we get, for the plateau
and the well
\begin{equation}
R_{0,L}^{x}=\frac{R_{x}^{0}\left(R_{L}^{0}-R_{x}^{0}\right)}{R_{L}^{0}}\asymp\mathcal{Z}\times\begin{cases}
\min\left\{ x,\, L+L^{c-b}-x\right\}  & \text{ for }x\le L-2\\
L^{c-b} & \text{ for }x=L-1
\end{cases}\label{eq:parall}
\end{equation}
The maximal local time in the plateau is where the resistance $R_{0,L}^{x}$,
the parallel between the resistances $R_{0}^{x}$and $R_{L}^{x}$,
is maximal.
\begin{itemize}
\item if $1>c-b$, the resistance of the plateau is larger than that of
the well and of the peak. Depending on the depth of the well, the
maximal resistance is acheived either in the middle of the plateau
or in the well. Indeed, an upper bound for $\max_{1\le x\le L-3}R_{x,L}^{0}$is
obtained by maximizing the r.h.s. of the first equality in (\ref{eq:parall});
as a function of $R_{x}^{0}$, this quantity has a maximum for $R_{x}^{0}=R_{L}^{0}/2$;
a lower bound for $\max_{1\le x\le L-3}R_{x,L}^{0}$is $R_{0,L}^{\lceil L/2\rceil}$.
Thus,
\[
\max_{1\le x\le L-3}\mathbb{E}(\xi_{0,L}^{x})\asymp L.
\]
The well $x=L-2$ has a large invariant measure that may compensate
for the small resistance. Its local time is $\mathbb{E}(\xi_{0,L}^{L-2})=\mu(L-2)\, R_{0,L}^{L-2}\asymp L^{c}$.
The local time of the peak $x=L-1$ is always negligible: $\mathbb{E}(\xi_{0,L}^{L-1})=\mu(L-1)R_{0,L}^{L-1}\asymp1$.
Therefore, the maximum local time is either in the middle of the plateau
or in the well and
\[
\max_{1\le x\le L-1}\mathbb{E}(\xi_{0,L}^{x})\asymp L+L^{c}.
\]

\item if $1<c-b$, the resistance of the plateau is negligible and the well
always wins:
\[
R_{0,L}^{x}\asymp\begin{cases}
x & \text{ for }x\le L-3\\
L^{c-b} & \text{ for }x\ge L-2
\end{cases}
\]
 Thus,
\end{itemize}
\[
\max_{1\le x\le L-3}\mathbb{E}(\xi_{0,L}^{x})\asymp\mu(L-2)R_{0,L}^{L-2}.
\]

Altogether,

\begin{equation}
\max_{x\le L-1}\mathbb{E}(\xi_{0,L}^{x})=\max_{x\le L-1}\mu(x)\, R_{0,L}^{x}\asymp\begin{cases}
\mathbb{E}(\xi_{0,L}^{\lceil\frac{L}{2}\rceil})\asymp L & \text{ for }c<1\\
\mathbb{E}(\xi_{0,L}^{L-2})\asymp L^{c} & \text{ for }1\le c\le b+1\\
\mathbb{E}(\xi_{0,L}^{L-2})\asymp L^{b+1} & \text{ for }c>b+1
\end{cases}.\label{eq:stimaLT}
\end{equation}

\subsubsection{Hitting times}


In the one-dimensional case, $\mathbb{P}\left(\tilde \tau_{k}^{0}<\tilde \tau_{L}^{0}\right)=1$
for all $k\in[1,L-1]$. By (\ref{eq:Txi}), (\ref{eq:r}), we see
that

\begin{gather}
\mathbb{E}(\tau_{L}^{0})=\sum_{k=0}^{L-1}\mathbb{E}\left(\xi_{L}^{k}\right)\mathbb{P}\left(\tilde \tau_{k}^{0}<\tilde \tau_{L}^{0}\right)=\sum_{k=0}^{L-1}\mu(k)R_{L}^{k}\nonumber \\
\mathbb{=E}(\xi_{L}^{0})+\sum_{k=1}^{L-3}\left(2(L-2-k)+2L^{c-b}\right)+\pi(L-2)4L^{c-b}+\pi(L-1)2L^{c-b}
\end{gather}

by (\ref{eq:LT}), (\ref{eq:pi})
\begin{align}
\mathbb{E}(\tau_{L}^{0})\asymp & \left(L^{a}\left(L+L^{c-b}\right)+L^{2}+L^{c-b+1}+L^{c}+1\right)(1+o(1))=\nonumber \\
 & \left(L^{a}+L\right)\left(L^{c-b}+L\right)(1+o(1)),\label{eq:T-1}
\end{align}
 where we used $a>b.$

The computation of $\mathbb{E}(\tau_{0,L}^{x})$ is slightly more
intricate:

In one dimension, $\mathbb{P}\left(\tilde \tau_{k}^{x}<\tilde \tau_{0,L}^{x}\right)$
is equal to $\mathbb{P}\left(\tilde \tau_{k}^{x}<\tilde \tau_{L}^{x}\right)$ if
$k\le x$ and to $\mathbb{P}\left(\tilde \tau_{k}^{x}<\tilde \tau_{0}^{x}\right)$
if $k>x$.
\[
\mathbb{P}\left(\tilde \tau_{k}^{x}<\tilde \tau_{0,L}^{x}\right)=\begin{cases}
\mathbb{P}\left(\tilde \tau_{k}^{x}<\tilde \tau_{L}^{x}\right)=\frac{R_{L}^{k}}{R_{L}^{x}} & \text{ if }k\le x\\
\mathbb{P}\left(\tilde \tau_{k}^{x}<\tilde \tau_{L}^{x}\right)=\frac{R_{0}^{k}}{R_{L}^{x}} &  \text{ if }k>x
\end{cases}
\]

By (\ref{eq:Txi}), for $x\in[1,\, L-1]$,
\begin{eqnarray}\label{tx0L}
\mathbb{E}(\tau_{0,L}^{x}) & =&\sum_{k=1}^{x-1}\mathbb{E}(\xi_{0,L}^{k})\mathbb{P}\left(\tilde \tau_{k}^{x}<\tilde \tau_{L}^{x}\right)+\mathbb{E}(\xi_{0,L}^{x})+\sum_{k=x+1}^{L-1}\mathbb{E}(\xi_{0,L}^{k})\mathbb{P}\left(\tilde \tau_{k}^{x}<\tilde \tau_{0}^{x}\right)\nonumber \\
 & =&\displaystyle
 \frac{R^x_{L}\sum_{k=1}^{x-1}\mu(k)R^k_{0}+
 \mu(x)R^0_{x}R^x_{L}+R^0_{x}\sum_{k=x+1}^{L-1}
 \mu(k)R^k_{L}}{R^0_{L}}\label{eq:Tn0L}
\end{eqnarray}

\begin{itemize}

\item Let us consider first the plateau $x\le L-3.$ By (\ref{eq:R0n})
\begin{equation}
\sum_{k=1}^{x-1}\mu(k)R^k_{0}\asymp x^{2},\label{eq:s1}
\end{equation}
and by (\ref{eq:RnL})
\begin{eqnarray}
\sum_{k=x+1}^{L-1}\mu(k)R^k_{L} & = & \sum_{k=x+1}^{L-3}\left(2(L-2-k)+4L^{c-b}\right)+2L^{c}+2\nonumber \\
 & \asymp & (L-x)^{2}+L^{c-b+1}+L^{c}.\label{eq:s2}
\end{eqnarray}
Plugging (\ref{eq:R0n},\ref{eq:RnL},\ref{eq:s1},\ref{eq:s2}) into
(\ref{eq:Tn0L}), we get
\begin{equation}
\mathbb{E}(\tau_{0,L}^{x})\asymp x\frac{\left(L-x+L^{c-b}\right)x++(L-x)^{2}+L^{c-b+1}}{L+L^{c-b}},\label{eq:Tlike}
\end{equation}
where we used $\mu(L-2)R_{L}^{L-2}\asymp L^{b}L^{c-b}$.
A little algebra
shows that
\begin{flalign*}
\mathbb{E}(\tau_{0,L}^{x})\asymp &  x L-\frac{L^{c-b}-L}{L^{c-b}+L}x^{2}\asymp\begin{cases}
xL-x^{2} & \text{ for }1\ge c-b\\
xL+x^{2} & \text{ for }1<c-b
\end{cases}
\end{flalign*}

In both cases,
\[
\max_{1\le x\le L-3}\mathbb{E}(\tau_{0,L}^{x})\asymp L^{2}
\]
\item if $x=L-2$, the three terms in the numerator of
l.s.h. of \eqref{tx0L} become:
\be{term1}
R^{L-2}_{L}\sum_{k=1}^{L-3}\mu(k)R^k_{0}
\asymp
\mathcal{Z} L^{2+c-b}
\ee
\be{term2}
\mu(L-2)R^0_{L-2}R^{L-2}_{L}
\asymp
\mathcal{Z} L^b L L^{c-b}
=
\mathcal{Z} L^{c+1}
\ee
\be{term3}
R^{0}_{L-2}\mu(L-1)R^{L-1}_{L}
\asymp
\mathcal{Z} L^{b-c} L L^{c-b}
=
\mathcal{Z} L,
\ee
where we used \ref{eq:RnL}.

By using \eqref{tx0L} and \ref{eq:RnL}, we estimate
\begin{flalign} \label{eq:Tlike2}
\mathbb{E}(\tau_{0,L}^{L-2})\asymp &  \frac{L^{c-b+2}+L^{c+1}}{L^{c-b}+L}x^{2}\asymp\begin{cases}
L^{c-b+1}+L^{c} & \text{ for }1\ge c-b\\
L^2+L^{b+1} & \text{ for }1<c-b
\end{cases}
\end{flalign}
\item if $x=L-1$, it is easy to see that
$$
\mathbb{E}(\tau_{0,L}^{L-1})
\le
\mathbb{E}(\tau_{0,L}^{L-2}).
$$
\end{itemize}
Putting together the three cases, we get
\begin{flalign} \label{DAFARE}
\max_{x \in \{1,L-1\}} \mathbb{E}(\tau_{0,L}^{x})\asymp
\begin{cases}
L^{b+1} & \text{ for } b>1 \text{ and } c > b+1 \\
L^2 & \text{ otherwise. }
\end{cases}
\end{flalign}




\end{document}